\input psfig.sty
\input psfig.sty

\tolerance=4000
\magnification=\magstep1
\parskip=4pt plus 6pt

\font\eightrm=cmr8
\font\bigbf=cmbx12

\font\vbigbf=cmbx12 scaled \magstep1

% \font\tencmr=cmr10
% \font\eightcmr=cmr8
% \font\sevencmr=cmr7
% \font\fivecmr=cmr5
% \newfam\cmrfam
% \textfont\cmrfam=\tencmr
% \scriptfont\cmrfam=\sevencmr
% \scriptscriptfont\cmrfam=\fivecmr
% \def\cmr{\fam\cmrfam}
% 
% 

\font\eusml=eusm10
\font\eusms=eusm7
\font\eusmf=eusm5
\newfam\eusmfam
\textfont\eusmfam=\eusml
\scriptfont\eusmfam=\eusms
\scriptscriptfont\eusmfam=\eusmf
\def\eusm{\fam\eusmfam}
% 
% \font\tenbf=cmbx10
% \font\eightbf=cmbx8
% \font\sevenbf=cmbx7
% \font\fivebf=cmbx5
% \newfam\cmbxfam
% \textfont\cmbxfam=\tenbf
% \scriptfont\cmbxfam=\sevenbf
% \scriptscriptfont\cmbxfam=\fivebf
% \def\cmbx{\fam\cmbxfam}
% 
% \font\tenmsbm=msbm10
% \font\sevenmsbm=msbm7
% \font\fivemsbm=msbm5
% \newfam\msbmfam
% \textfont\msbmfam=\tenmsbm
% \scriptfont\msbmfam=\sevenmsbm
% \scriptscriptfont\msbmfam=\fivemsbm
% \def\msbm{\fam\msbmfam}

\def\cmr{\rm}
\def\cmbx{\bf}

\def\nobf{\noindent \bf}

\def\meti#1{\par\indent\llap{#1\enspace}\ignorespaces}

%\hfill\vrule height3pt width6pt depth3pt\smallskip}}

\def\itembu#1{\item{$\circ$}{#1}}

\def\im{{\cmbx im}}

\def\xx{{\cmbx x}}

\def\Fa{{{\cal F}_a}}
\def\Fb{{{\cal F}_b}}
\def\Fc{{{\cal F}_c}}
\def\Fx{{\cal F}}
\def\Cx{{\cal C}}

\def\Fbd{{\bf F}}
\def\Sbd{{\bf S}}
\def\Fabd{{{\eusm F}_a}}
\def\Fbbd{{{\eusm F}_b}}
\def\Fcbd{{{\eusm F}_c}}

\def\wz{{\cal W}}

\def\sgn{{\cmr sgn}}

\def\emptyset{\varnothing}

\mathchardef\braceld="37A
\mathchardef\bracerd="37B

\def\m@th{\mathsurround=0pt }

\def\downbrakfill{$\m@th
\;\braceld\leaders\vrule\hfill\bracerd\;$}

\def\overbrak#1{\mathop{\vbox{\ialign{##\crcr\noalign{\kern3pt}
\downbrakfill\crcr\noalign{\kern3pt\nointerlineskip}
$\hfil\displaystyle{#1}\hfil$\crcr}}}\limits}

\def\overline{\bar}

\def\arctanh{\hbox{arctanh}}
\def\aux{\hbox{\fiverm aux}}
\def\fancyone{{\bf 1}}
\input amssym.def
\input amssym.tex

\def\BB{{\Bbb B}}

\def\RR{{\Bbb R}}
\def\Ss{{\Bbb S}}
\def\SS{{\Bbb S}}

\def\emptyset{\varnothing}

\def\fullsquare{{\vrule height6pt width6pt depth0pt}}
\def\emptysquare{{\hbox{\vrule height6pt width0.6pt depth0pt%
\vbox{\hrule height0.6pt width4.8pt depth0pt%
\vglue4.8pt%
\hrule height0.6pt width4.8pt depth0pt}%
\vrule height6pt width0.6pt depth0pt}}}

\def\qed{\unskip\nobreak
\hfil\penalty50\hskip1.75em\null\nobreak\hfil\fullsquare
{\parfillskip=0pt \finalhyphendemerits=0 \par}}
\def\qede{\unskip\nobreak
\hfil\penalty50\hskip1.75em\null\nobreak\hfil\emptysquare
{\parfillskip=0pt \finalhyphendemerits=0 \par}}

% From "testfont.tex",
% a testbed for font evaluation (see The METAFONTbook, Appendix H)
%
% Translated, isolated and expanded by Nicolau C. Saldanha

\newcount\m \newcount\n \newcount\p \newdimen\dim

\def\today{\ifcase\month\or
  January\or February\or March\or April\or May\or June\or
  July\or August\or September\or October\or November\or December\fi
  \space\number\day, \number\year}

\def\hoje{\number\day\ de
  \ifcase\month\or
  janeiro\or fevereiro\or mar{\c c}o\or abril\or maio\or junho\or
  julho\or agosto\or setembro\or outubro\or novembro\or dezembro\fi
  \ de \number\year}

\def\japtoday{\number\year.\twodigits\month.\twodigits\day}

\def\hours{\n=\time \divide\n 60
  \m=-\n \multiply\m 60 \advance\m \time
  \twodigits\n:\twodigits\m}

\def\twodigits#1{\ifnum #1<10 0\fi \number#1}

\def\datepages{\footline={{\hss\tenrm\folio\hss\hbox to 0pt{
\hidewidth\fiverm\japtoday}}}}

\def\hu{{H^1}}

\def\ints{\int}

\datepages

\centerline{\vbigbf Morin singularities and global geometry}
\vskip4pt
\centerline{\vbigbf in a class of ordinary differential operators}
\bigskip
\centerline{\it Iaci Malta, Nicolau C. Saldanha and Carlos Tomei
\footnote{}{\eightrm Research supported by MCT and CNPq, Brazil}}

\bigskip \bigskip

{

\narrower

{\nobf Abstract: }
We consider the operator $F(u) = u' + f(t,u(t))$ acting on periodic
real valued functions.
Generically, critical points of $F$ are 
infinite dimensional Morin-like singularities
and we provide operational characterizations
of the singularities of different orders.
A global Lyapunov-Schmidt decomposition of $F$
converts $F$ into adapted coordinates,
$\Fbd(\tilde v, \overline u) = (\tilde v, \overline v)$,
where $\tilde v$ is a function of average zero
and both $\overline u$ and $\overline v$ are numbers.
Thus, global geometric aspects of $F$ reduce to the study
of a family of one-dimensional maps:
we use this approach to obtain normal forms for several nonlinearities $f$.
For example, we characterize autonomous nonlinearities
giving rise to global folds and, in general,
we show that $F$ is a global fold if all critical points are folds.
Also, $f(t,x) = x^3 - x$, or, more generally,
the Cafagna-Donati nonlinearity, yield global cusps;
for $F$ interpreted as a map between appropriate Hilbert spaces,
the requested changes of variable to bring $F$ to normal form
can be taken to be diffeomorphisms.
A key ingredient in the argument is the contractibility of both
the critical set and the set of non-folds for a generic autonomous
nonlinearity.
We also obtain a numerical example of a polynomial $f$ of degree 4
for which $F$ contains butterflies (Morin singularities of order 4)---%
it then follows that $F(u) = v$ has six solutions for some $v$.



}

\bigbreak

{\noindent\bigbf Introduction}

\nobreak\medskip\nobreak

In this paper we consider the differential equation
$$u'(t) + f(t,u(t)) = g(t), \eqno{(\ast)}$$
where the unknown $u$ is a real function on $\SS^1$ and the 
nonlinearity $f:\SS^1 \times \RR \to \RR$ can assume a
number of forms.

Our approach is to study the global geometry of the operator
$$\eqalign{F: B^1 &\to B^0 \cr u &\mapsto u' + f(t,u)}$$
where the domain is either $C^1(\SS^1)$ (the Banach  space of 
periodic functions with continuous derivatives) or the Hilbert
space $H^1(\SS^1)$ of periodic functions with square integrable
derivative. Ideally, we search for global changes of variables
in both domain and image taking the operator $F$ to a simple normal form.
This goal has been achieved in previous occasions,
starting with the seminal work of Ambrosetti and Prodi ([AP])
and its geometric interpretation by  Berger and Church ([BC]),
who showed that the operator associated to a certain nonlinear
Dirichlet problem gives rise to a global fold between infinite dimensional
spaces. 
Topological global cusps have appeared already in
operators related to partial differential equations with a parameter
([BCT], [CDT]).
Closer to the subject of this paper, McKean and Scovel ([McKS]) 
showed that the operator $F$
for $f(t, x) = x^2$ (or more generally, for convex nonlinearities)
is also a global fold, and raised the question of the global nature of $F$
for $f(t, x) = x^3 - x$.
The same question was asked by Cafagna and Donati ([CD])
and Church and Timourian ([CT]), who state and prove
partial results for the more general Cafagna-Donati equation ([CD]),
for which $f(t,x) = ax + bx^2 + cx^{2k+1}$
for appropriate choices of $a$, $b$ and $c$.
In Theorem 5.1 and Corollary 5.5,
we show that these nonlinearities indeed obtain global cusps.
With some additional effort, we show that in Hilbert spaces
the requested global changes of variables can be taken to be smooth.

Actually, the operator $F$ is simple enough that substantial insight
into its global geometry can be obtained with the rather mild hypothesis
of $f$ being {\it tame} (see definition before Theorem 1.2).
First, we construct a global Lyapunov-Schmidt type decomposition of $F$.
Split a function
$u$ as a sum of a function of average zero ($\tilde u$) and
a constant ($\overline u$) and decompose domain and 
image accordingly: $B^i = \tilde B^i \oplus \overline B^i$.
Writing the action of $F$ as
$$u = (\tilde u, \overline u) \mapsto v = (\tilde v, \overline v),$$
we show in Theorem 1.2 that, for each $\overline u$, the correspondence
$\tilde u \mapsto \tilde v$ is a global diffeomorphism.
This provides a change of coordinates in the domain of $F$
bringing it to (global) {\sl adapted coordinates}
$$(\tilde v , \overline u) \mapsto (\tilde v , \overline v).$$
We immediately obtain that the inverse images of vertical lines
under $F$ are {\sl fibres}, curves foliating the domain and
intersecting every horizontal plane exactly once and transversally.
The study of $F$ in a sense boils down to the study of its behaviour
on the fibres: for example, $f(t, x) = x^2$ produces a fold
on every fibre and thus $F$ is a global fold. 

Only the tameness hypothesis on the behaviour of $f$ at infinity
is necessary to obtain adapted coordinates:
in particular, we obtain some results about the
global geometry of $F$ even when it is not proper.
From Proposition 1.4, properness of $f$ implies properness of $F$
but, from Proposition 4.1, the converse is false.
Adapted coordinates combined with properness make clear the possibility
of definining a {\sl topological degree} for $F$:
the degree of $F$ is just the degree of any of its restrictions to fibres.

It is easy to see that (generically) $S_1$,
the critical set of $F$, is a manifold.
Rather surprisingly, the global geometry of $S_1$
does not depend on the nonlinearity:
generically, it is connected and contractible (Corollary 1.9).
This follows from a more general theorem (Theorem 1.8 or [MST])
on the contractibility of regular level sets of a class
of functionals defined by integration.
From contractibility, by topological arguments
often using the infinite dimension of the spaces involved ([Ka], [Ku], [S]),
there is a change of variables in the domain of $F$
taking $S_1$ to a closed hyperplane;
in the Hilbert case, this change of variables
can be taken to be a diffeomorphism but
in the Banach case, it is merely a homeomorphism.

We then proceed to study the critical points of $F$ in detail.
From adapted coordinates, $\ker DF$ has dimension 1 at critical points
and $\im DF$ is then a closed subspace of codimension 1.
This restricts considerably the possible nature of
a generic critical point of $F$:
it has to be an infinite dimensional {\sl Morin singularity} ([M]).
More precisely,
after changes of coordinates $F$ near a generic singularity $u$ of 
can be written as
$$(Z,x_1,\ldots,x_{k-1},y) \mapsto
(Z,x_1,\ldots,x_{k-1},y^{k+1} + \sum_{i = 1, \ldots, k-1} x_i y^i)$$
near zero, where $Z$ is an element of an infinite dimensional space and
$x_i$ and $y$ are real numbers.
The integer $k$ is the {\sl order} of the singularity:
folds and cusps are Morin singularities of orders 1 and 2.
Morin's classification and proof carry over to the infinite dimensional case by
making use of a version of the Malgrange preparation theorem with an
(infinite dimensional) parameter: this approach has been used in [CDT] to
obtain a characterization of infinite dimensional cusps.
The description of a Morin singularity is 
given more explicitly in Propositions 2.1 and 2.2 
in terms of a collection of functionals $\Sigma_i, i=1,\ldots$: 
at a singular point of order $k$, the first $k$ functionals have to be zero 
and some transversality relations have to hold. 
Given $v$, we may define a {\sl return map} $\rho_v$ 
taking $x_0$ to $x_1$ if a (possibly non-periodic) solution $u$ of $(\ast)$
satisfies $u(0) = x_0$, $u(1) = x_1$.
In Proposition 2.3, we relate the order of a singularity $u$
to the order of contact between $\rho_{F(u)}$ and the identity at $u(0)$.

In the autonomous case, when the nonlinearity does not depend on $t$,
also $S_2$, the set of critical points which
are not folds, is (generically) a connected contractible manifold. 
To show this, we need again Theorem 1.8 and
an Lemma 3.1, stating that $S_1$ and $S_2$
are diffeomorphic to the simpler sets $\hat S_1$ and $\hat S_2$,
critical and non-fold points of the {\sl simplified operator}
$$\eqalign{
\hat F:B^1 &\to B^0. \cr
u &\mapsto u' + \ints f(t,u(t)) dt\cr}$$
Again, contractibility yields a change of variables
flattening $S_1$ and $S_2$.

The functionals $\Sigma_i$ are rather complicated and we
do not know of a simple procedure to decide if singularities of a given
order exist for a fixed $f$.
However, in the autonomous case,
we describe in Lemma 3.5 a necessary (and essentially sufficient) criterion
for the existence of singularities of order $k$
for the simplified operator $\hat F$.
In the same lemma, we show that if $\hat F$
has a singularity of order $k$ then $F$ also does.

In Section 4, we consider some special types of functions $f$:
if $f$ is either monotonic or convex 
for each value of the first coordinate $t$,
we give a global description of the behaviour of $F$.
Even though some of the results are simple or well known (from [McKS]),
they provide a convenient introduction to our approach
of studying $F$ fibre by fibre.
Using Lemma 3.5, we give a criterion (Theorem 4.4)
for autonomous nonlinearities
to decide whether the operator $F$ is a global fold.
In particular (Corollary 4.5), polynomial non-convex nonlinearities $f$
give rise to operators $F$ with $S_2 \ne \emptyset$
but there are fast-growing non-convex nonlinearities
for which the operator is a global fold.
Also, local behaviour characterizes global folds (Theorem 4.6):
generically, if all singularities of $F$ are folds then 
$F$ is a global fold.

In Section 5, the autonomous nonlinearities $f$ satisfy $f''' \ge 0$
with isolated zeros.
The related operator $F$ is then a global cusp:
here we make full use of our techniques.
In this case there does not seem to be an explicit description 
of the requested (global) changes of variables:
their existence follows by topological arguments similar
to those used in the study of the sets $S_1$ and $S_2$.  
Lemma 5.4 is a global parametrized version of Whitney's normal form
for cusps ([W]); the proof appears to be cumbersome but its main
difficulty lies in verifying that Whitney's construction
can be performed smoothly in a parameter.
We present only a sketch of argument
and we thank John Mather for helpful discussions.

We finish the paper with an example of a different kind.
The results in sections 4 and 5 are enough to show that if $f(t,x)$
is a polynomial in $x$ of degree $d \le 3$
(with coefficients depending on $t$ and non-zero coefficient of highest degree),
then the related operator $F$ is a diffeomorphism,
a global fold or a global cusp.
In this case, thus, equation $(\ast)$ has at most $d$ periodic solutions.
The number of solutions of $(\ast)$ when $f$ is such a polynomial
was considered by Pugh, Lins Neto and Smale ([L])
who proved the bounds above and that
the number of solutions may be arbitrarily large for $d = 4$.
We instead exhibit a numerical example of an autonomous polynomial
$f$ of degree four and a function $u$ which is a Morin
singularity of order four (a {\sl butterfly}). This was accomplished by
requesting that $u$ be a root of the first four functionals $\Sigma_i$.
By the normal form of $F$ at a butterfly, there are points $g$ near $F(u)$
with five pre-images; one is presented. By a degree-theoretic 
argument, such a (regular) point ought to have an even number of
pre-images, and we verified by solving the differential equation
with a Runge-Kutta method that there are exactly six initial conditions
giving rise to periodic solutions.

\bigbreak

\bigbreak

{\noindent\bigbf 1. Adapted coordinates and the critical set}

\nobreak\smallskip\nobreak

We consider the smooth nonlinear operator
$F: B^1 \to B^0$ given by $$F(u)(t) = u'(t) + f(t,u(t)),$$
where $f: \SS^1 \times \RR \to \RR$ is a smooth function.
Here $B^1$ and $B^0$ can be chosen in two different ways.
In the $H$ case, they are the Sobolev spaces
$B^1 = H^1 = H^1(\SS^1;\RR)$
(the periodic absolutely continuous real valued functions
with derivative in $L^2$) and $B^0 = H^0 = L^2(\SS^1;\RR)$.
In the $C$ case,
$B^1 = C^1 = C^1(\SS^1;\RR)$ and $B^0 = C^0 = C^0(\SS^1;\RR)$.
For notational convenience,
inner products are to be interpreted in the $L^2$ sense
even in other spaces.
An interesting special situation is the {\sl autonomous} case,
in which $f$ does not depend on the $t$ coordinate.
We denote the partial derivative of $f$ with respect to
the second variable by $D_2f$.

Proposition 1.1 below obtains a formula for $DF$ at arbitrary points,
a description of the critical set $S_1F$ and
a Lyapunov-Schmidt decomposition for the operator $F$
in a neighbourhood of a critical point in the domain.
In Theorem 1.2 we show the existence
of a convenient global decomposition of $F$.

Recall the familiar Green kernel $k(x) = x - \lfloor x \rfloor - 1/2$
(where $\lfloor x \rfloor$, following Knuth,
is the largest integer not larger than $x$).
If $h$ is periodic (with period 1)
then $h_1(t) = \ints k(s-t) h(s) ds$ is also periodic and
$h_1'(t) = h(t) - \ints h(s) ds$, a function of average 0,
so that $k$ is a kernel for the inverse,
restricted to functions of average zero, of the derivative.

{\nobf Proposition 1.1: }
{\sl 
The derivative
$$(DF(u) v)(t) =  v'(t) + D_2f(t,u(t)) v(t)$$
is a Fredholm operator of index 0 from $B^1$ to $B^0$;
furthermore, $\ints D_2f(s,u(s)) ds$ is the unique real eigenvalue
of $DF(u)$, which is simple,
with corresponding eigenvector is
$$w_u(t) = e^{-\ints k(s-t) D_2f(s, u(s)) ds}.$$
In particular, the critical set of $F$ is
$$S_1F = \left\{{u \in B^1 | \ints D_2f(t,u(t)) dt = 0}\right\}.$$
The subspace $\langle 1/w \rangle^\perp$ has codimension 1,
is transversal to $\langle w \rangle$ and
is also invariant under $DF(u)$.
Thus, the restriction 
$$DF(u): \langle 1/w \rangle^\perp \subset B^1 \to
\langle 1/w \rangle^\perp \subset B^0$$
is bijective.
}

By an eigenvector of $DF(u): B^1 \to B^0$ we mean a solution
of $DF(u)v = v' + D_2f(t,u(t)) v = \lambda v$;
by standard regularity arguments,
solutions of this equation are always in $B^1$.

{\nobf Proof: }
The formula for the derivative is straightforward.
The expression for $w$ follows from the explicit solution of
the first order periodic linear ODE and $1/w$
is the only real eigenvector of the (adjoint) operator
$$v \mapsto - v' + D_2f(t, u(t)) v.$$
\qed

Let $B^1 = \tilde B^1 \oplus \langle 1 \rangle$ and
$B^0 = \tilde B^0 \oplus \langle 1 \rangle$
(the tilde denotes integral equal to 0)
defining complementary projections 
$\Pi_{\tilde B}$ and $\Pi_{\overline B}$.
More concretely,
$\tilde u = \Pi_{\tilde B} u = (u - \int u) + \int u$ and
$\overline u = \Pi_{\overline B} u = \int u$
(we omit the domain of integration when it is $\SS^1$).
Notice that $\langle 1 \rangle$ is always transversal
to $\langle 1/w \rangle^\perp$
and that $\langle w \rangle$ is likewise transversal to
$\tilde B^i$.

A function $f: \Ss^1 \times \RR \to \RR$ is
{\it wild at} $+\infty$ (resp. $-\infty$) if
$$ \int_I {{ds} \over {\max(1,\sup_{t \in \Ss^1} f(t,s))}} < +\infty, \quad
\int_I {{ds} \over {\max(1,\sup_{t \in \Ss^1} (- f(t,s)))}} < +\infty $$
for $I = [0,+\infty)$ (resp. $I = (-\infty, 0]$);
$f$ is {\it tame} if not wild at $\pm\infty$.
Loosely, $f$ being tame implies that a solution $u$
can not go very far and come back in bounded time.

{\nobf Theorem 1.2: }
{\sl Let $f: \Ss^1 \times \RR \to \RR$ be a tame nonlinearity.
Let $F(\tilde u + \overline u) = \tilde v + \overline v$.
The map $\Psi: B^1 \to B^0$,
$\Psi(u) = \tilde v + \overline u$, is a (global) diffeomorphism.}

The tameness hypothesis cannot be discarded.
For $f(t,u) = 2\pi \cos(2\pi t) \cosh^2(u)$,
there are no periodic functions $u$ such that
$u'(t) + f(t,u(t))$ is constant.
In particular, the point $0 \in B^0$
is not in the image of the map $\Psi$ above.
Indeed, the solutions of $u'(t) + f(t,u(t)) = 0$ are
$$ u = -\arctanh\left( \sin(2\pi t) + C \right), \quad C \in (-2,2). $$
For $C = 0$, consider the solutions $u_-$ and $u_+$
on disjoint domains $(-1/4,1/4)$ and $(1/4,3/4)$.
Notice that $u_-$ (resp. $u_+$) is strictly decreasing (resp. increasing)
with absolute value tending to infinity at the endpoints of the domain.

% \begin{figure}[ht]
% \psfrag{u}{$u$}
% \psfrag{t}{$t$}
% \psfrag{u0}{$u_\nu$}
% \psfrag{u+}{$u_+$}
% \psfrag{u-}{$u_-$}
% %\begin{center}
\smallskip
\centerline{\psfig{height=50mm,file=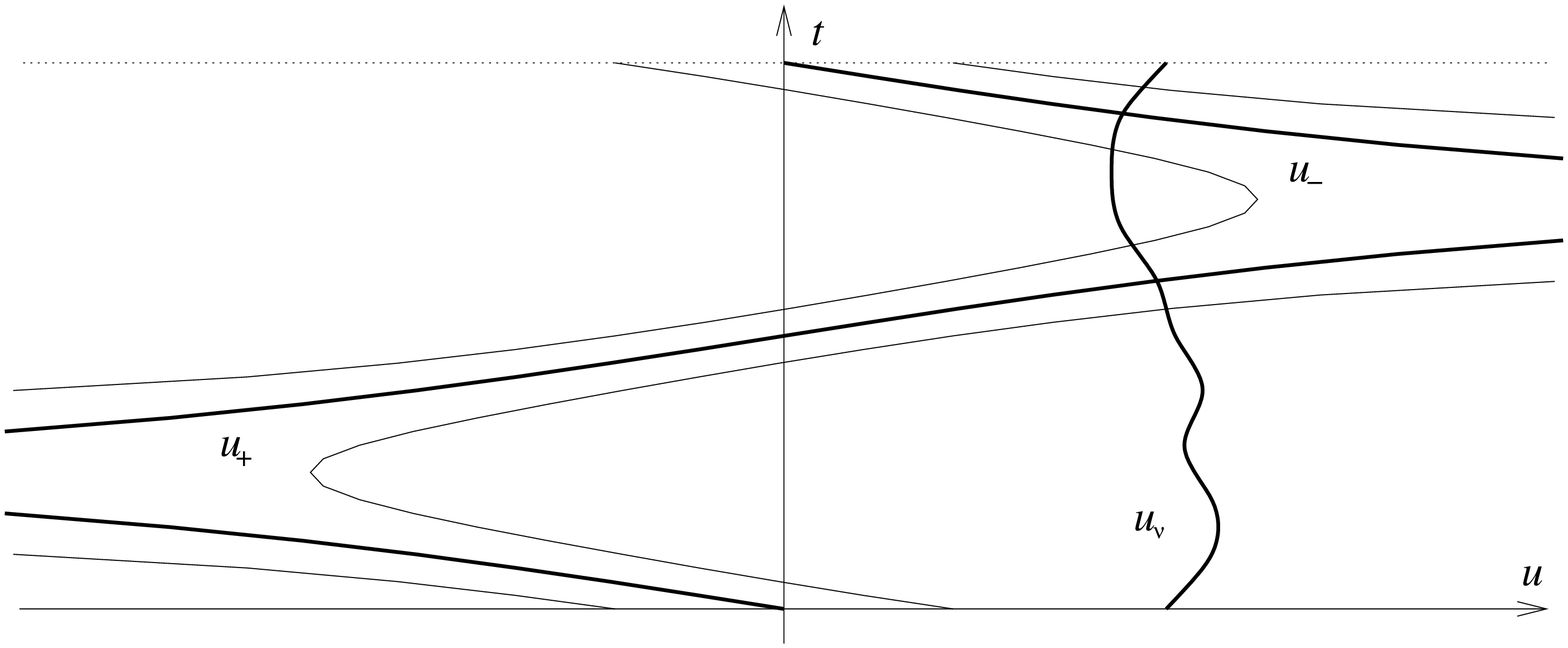}}
% %\end{center}
% \caption{Solutions for $C = 0$ and $C = \pm 0.3$}
% \label{fig:curves}
% \end{figure}

The graph of any periodic function $u_\nu$ must cross the graphs
of both $u_-$ and $u_+$ at times $t_-$ and $t_+$, respectively,
for which $u_\nu'(t_-) + f(t_-,u_\nu(t_-)) \ge 0$ and 
$u_\nu'(t_+) + f(t_+,u_\nu(t_+)) \le 0$.
If $u_\nu'(t) + f(t,u_\nu(t)) = \nu$ for all $t$ then,
from the conditions above, $\nu = 0$.
This, however, implies that
$u_\nu$ must equal both $u_-$ and $u_+$, a contradiction.

\smallskip

We provide some equivalent, more geometric, readings
for the rather dry statement of Theorem 1.2.
The following diagram may be helpful:
\smallskip
\centerline{\psfig{width=55mm,file=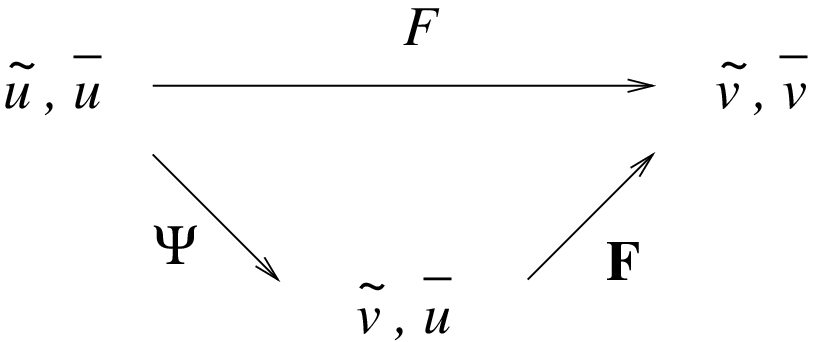}}

The map $\Fbd = F \circ \Psi^{-1}: B^0 \to B^0$
takes $(\tilde v, \overline u)$ to
$(\tilde v, \overline v) = (\tilde v, \phi(\tilde v,\overline u))$.
Horizontal hyperplanes $\tilde B^1 + \{c\}$
are injectively taken by $F$ onto {\sl sheets},
i.e., hypersurfaces intersecting each vertical line
$\{\tilde v\} + \langle 1 \rangle$ transversally and exactly once.
Equivalently, the inverse images under $F$ of the vertical lines
$\{\tilde v\} + \langle 1 \rangle$ foliate $B^1$ by {\sl fibres},
i.e., curves intersecting each horizontal hyperplane transversally
and exactly once;
we denote by $\tau_u$ the fibre containing $u$.
Let $T$ be the set of fibres:
from transversality, we may identify $T$ with any horizontal hyperplane
in the domain, in particular with $\tilde B^1$.
The set of vertical lines in the image is naturally identified
with $\tilde B^0$ and $F$ induces a diffeomorphism from $T$ to $\tilde B^0$.

The following proof applies to both cases
but certain complications are relevant only for the $H$ case.

{\nobf Proof: }
In order to invert a vertical line and obtain a fibre,
we consider the differential equation
$$u'(t) + f(t,u(t)) = \tilde v(t) + \nu,\eqno{(\ast)}$$
where $\tilde v \in \tilde B^0$ is fixed
and $\nu \in \RR$ is a parameter.
Local existence, uniqueness and continuous dependence on parameters
hold even when $\tilde v$ is only $L^2$.
Also, solutions cease to exist only by going to infinity.

{\sl Given $t_0$ and $u(t_0)$ there are $\epsilon > 0$, $\nu_+$ and $\nu_-$
such that the two solutions $u_+$ and $u_-$ of $(\ast)$ with
initial condition $u(t_0)$ satisfy:
\item{\rm ($+$)}{$u_+$ either goes to $+\infty$ at some time $t$,
$t_0 < t \le t_0 + \epsilon$ or satisfies $u_+(t_0 + \epsilon) > u_+(t_0)$.}
\item{\rm ($-$)}{$u_-$ either goes to $-\infty$ at some time $t$,
$t_0 - \epsilon \le t < t_0$ or satisfies $u_-(t_0 - \epsilon) < u_-(t_0)$.}
}

We discuss only ($+$): item ($-$) is analogous.
Notice that the claim is trivial in the $C$ case:
choose the parameter $\nu_+$ so that the derivative at time $t_0$
of $u_+$ is positive.
Clearly, if ($+$) is satisfied by some $\nu_+$,
it is satisfied by sufficiently positive $\nu_+$.

Solve $(\ast)$ for $\nu = 0$ to obtain a solution $u_0$
defined on $[t_0, t_0 + \epsilon]$.
Without loss, $u_0(t_0) > u_0(t_0 + \epsilon)$.
Let $u_{\aux}(t) =
u_0(t) + \epsilon^{-1}(u_0(t_0) - u_0(t_0 + \epsilon))(t - t_0)$.
We choose $\nu_+$ such that the vector field
$(1,- f(t,u) + \tilde v(t) + \nu_+)$ always crosses the graph of $u_{\aux}$
upwards, i.e.,
$-f(t,u_{\aux}(t)) + \tilde v(t) + \nu_+ > u_{\aux}'(t) =
-f(t,u_0(t)) + \tilde v(t) +
\epsilon^{-1}(u_0(t_0) - u_0(t_0 + \epsilon))$.
This is clearly possible since $f$ is continuous.
\qede

{\sl Given $u(0) = c$, there is some $\nu_+$ for which
the maximal solution $u_{\nu_+}$ of $(\ast)$
either goes to $+\infty$ at some time $t$,
$0 < t \le 1$, or satisfies $u_{\nu_+}(1) > c$.}
This follows from the previous claim by a compactness argument.
\qede

{\sl Given $u(0) = c$, there is a unique $\nu$
for which $(\ast)$ admits a periodic solution.}
Consider the set $A^+$ (resp. $A^-$) of $\nu$'s such that
$u_\nu$ goes to $+\infty$ (resp. $-\infty$)
or satisfies $u_\nu(1) \ge c$ (resp. $u(1) \le c$).
From the previous claim (and the obvious counterpart),
both sets are non-empty.
Set $\nu_0 = \sup A^- = \inf A^+$:
if $\nu_0 \in A^- \cap A^+$ then $u_{\nu_0}$ is periodic;
uniqueness follows from the local behaviour of the solutions.
We show that $\nu_0 \notin A^\mp$
implies that $f$ is wild at $\pm\infty$.
We consider the case $\nu_0 \notin A^-$.

If $u_{\nu_0}$ is defined in $[0,1]$ with $u_{\nu_0}(1) > c$
then continuous dependence implies that some open neighborhood of $\nu_0$
is contained in $A^+$, a contradiction.
Define $t_c \in (0,1]$ by
$$ \lim_{t \to t_c} u_{\nu_0}(t) = +\infty. $$
If $t_c = 1$, continuous dependence again implies that $\nu_0$
is in the interior of $A^+$; we therefore have $t_c < 1$.
Thus, for every $M \in \RR$ there exists $\nu < \nu_0$
and $t_\nu > t_c$ such that $u_\nu(t_c) > M$, $u_\nu(t_\nu) < c$.
Set
$$ I_{\tilde v} =
\{t \in [0,t_c] \;|\; -f(t,u_\nu(t)) \le \tilde v(t) + \nu \}. $$
For $t \in I_{\tilde v}$ we have $u_\nu'(t) \le  2\tilde v(t) + 2\nu$
and the Lebesgue measure $\mu(u_\nu(I_v))$ is bounded above by
$2|\nu| + 2||\tilde v||_{L^1}$.
Define $h: [c,M] \to [0,t_c]$ by
$h(s) = \inf\{t \in [0,t_c] \;|\; u_\nu(t) = s\}$.
Even though $h$ may have discontinuities,
it is strictly increasing and then, for almost all $s$,
$h$ is differentiable with $h'(s) = 1/u_\nu'(h(s))$.
Let $J_f =  [c,M] \smallsetminus u_\nu(I_v)$:
for $s \in J_f$, we have $h(s) \notin I_{\tilde v}$
and $h'(s) \ge -1/(2f(h(s),s))$.
Thus,
$$ t_c \ge \mu(h(J_f)) \ge \int_{J_f} h'(s) ds
\ge {{1}\over{2}} \int_{J_f} {{ds}\over{\max(1,\sup_{t \in \Ss^1}(-f(t,s)))}} $$
and therefore
$$ \int_c^M {{ds}\over{\max(1,\sup_{t \in \Ss^1}(-f(t,s)))}}  \le
2t_c + 2|\nu| + 2||\tilde v||_{L^1}. $$
Since this estimate holds for arbitrarily large $M$,
$$ \int_0^{+\infty} {{ds}\over{\max(1,\sup_{t \in \Ss^1}(-f(t,s)))}} 
< +\infty. $$
A similar argument for the interval $[t_c,t_\nu]$ yields
$$ \int_0^{+\infty} {{ds}\over{\max(1,\sup_{t \in \Ss^1}(f(t,s)))}} 
< +\infty, $$
implying that $f$ is wild at $+\infty$.

\qede

{\sl Given $a$, there is a unique $\nu$
for which $(\ast)$ admits a periodic solution $u$
with $\int u = a$.}
Consider all periodic solutions of $(\ast)$ as curves
in $\SS^1 \times \RR$.
Again, from the local behaviour of solutions, the curves are disjoint.
By the previous claim, the union of all such curves
contains the line $\{0\} \times \RR$.
By a similar argument applied to other lines,
the union of the curves is $\SS^1 \times \RR$
and the curves form a continuous foliation of the cylinder by circles.
Notice that circles in the foliation correspond
to points in the fibre 
$F^{-1}(\{\tilde v\} + \langle 1\rangle)$.
The integrals of the solutions are strictly increasing
as a function of $u(0)$.
Also, there are solutions with arbitrarily large
(positive or negative) integrals, since the area between two curves
goes to infinity as the initial condition of one of the curves does.
\qede

At this point, we have that, given $\overline u$ and $\tilde v$,
there is a unique $\tilde u$ with
$F(\tilde u + \overline u) = \tilde v + \overline v$
(for some $\overline v$).
Thus, the function $\Pi_{\tilde B^0} \circ F$ is a bijection
from any hyperplane $\tilde B^1 + \{a\}$ to $\tilde B^0$.
This function is clearly smooth and, as discussed before the
statement of this theorem, its derivative is always invertible.
By the inverse function theorem, these bijections are diffeomorphisms.
\qed

\goodbreak

{\nobf Remarks: }
\meti{1.} Theorem 1.2 is the counterpart of the usual
global domain decomposition found in the study of the equation
$\Delta u = f(u)$, with Dirichlet boundary conditions
and special resonance hypothesis on $f$ (see [AP]).
There, the task is simplified by the use of self-adjoint
spectral theory.
In our case, the derivative, unlike the Laplacian, is skew-symmetric,
with purely imaginary spectrum containing 0,
and the nonlinearity interacts at most with 0.
\meti{2.} The hard part in an eventual functional analytic proof of Theorem 1.2
is the properness (and hence, from local behaviour, bijectivity) of
the function $\Pi_{\tilde H^0} \circ F:
\tilde H^1 + \{\upsilon\} \to \tilde H^0$;
the necessary estimates seem to be simple only in the $C$ case.

For later use, we state as a lemma
some consequences of the proof of Theorem 1.2.

{\nobf Lemma 1.3: }
{\sl 
Assume $f$ to be a tame nonlinearity.
Fibres are parametrized by average, i.e.,
the function
$$\tau_{u_0} \to \RR, \quad u \mapsto \ints u(t) dt$$
is a diffeomorphism.
Let $u_a$ be the element of average $a$ in $\tau_{u_0}$.
Then
$$ \lim_{a \to +\infty} \min_t u_a(t) = +\infty, \quad
\lim_{a \to -\infty} \max_t u_a(t) = -\infty.  $$
Given $t_0 \in \SS^1$, the function
$$\tau_{u_0} \to \RR, \quad u \mapsto u(t_0)$$
is also a diffeomorphism.
}

The study of the (global and local) geometry of $F$
thus reduces to the study of $\Fbd = F \circ \Psi^{-1}$
and therefore of $\phi: B^0 \to \RR$.
The diffeomorphism $\Psi$ is said to provide $F$
with {\sl adapted coordinates}, i.e.,
$\Fbd(\tilde v, \overline u) =
(\tilde v, \phi(\tilde v, \overline u))$.
This change of variables is convenient
to the classification of critical points of $F$,
as we shall see in the next section.
Notice that we do not have formulae for
$\Fbd$ or $\Psi$ and have to make do with $w$ and
$\Phi(u) = (\Pi_{\overline B_0} \circ F)(u) =
(\phi \circ \Psi)(u) = \ints f(t, u(t)) dt$
(a somewhat cumbersome formula for $\wz$ is given in Lemma 1.5).

{\nobf Proposition 1.4: }
{\sl If $f$ is proper then the operator $F: B^1 \to B^0$ is proper.}

Notice that $f: \SS^1 \times \RR \to \RR$ is proper
if and only if $|f(t,x)|$ goes to infinity when $(t,x)$ does.
In particular, if $f$ is proper, the restriction of $F$ to a fibre
takes infinity to infinity.
Also, if $f$ is proper then $f$ is necessarily tame.

{\nobf Proof: }
Any compact set $K \subseteq B^0$ is contained in the product of
its (compact) projections, so without loss, $K$ can be taken
to be the product of a compact set $\tilde K \subseteq \tilde B^0$
and an interval $[-k,k]$.
Since $F \circ \Psi^{-1}$ is of the form
$(\tilde v,\overline u) \mapsto (\tilde v, \phi(\tilde v,\overline u))$,
the compactness of the preimage of $K$
under $F \circ \Psi^{-1}$ (and hence under $F$)
follows from the boundedness of $\overline u$ in $F^{-1}(K)$
or the uniform boundedness of $u \in F^{-1}(K)$.

In the $C$ case, consider $u$ at its global extrema:
there, $v(t) = u'(t) + f(t,u(t)) = f(t,u(t))$
and properness of $f$ gives us the required uniform bound.
For the $H$ case,
assume by contradiction that there are $u_n \in F^{-1}(K)$ with
$u_n(t_n) > 2^n$, where $t_n$ is the global maximum of $u_n$, and,
without loss of generality,
that $f$ is positive for large positive $x$.
Consider the intervals $I_n = (t_n - 1/10, t_n)$.
If $u_n(t) > 2^{n-1}$ for all $t$ in $I_n$
for all sufficiently large $n$ then
$\int_{I_n}{v_n(t) dt} = u_n(t_n) - u_n(t_n - 1/10) +
\int_{I_n}{f(t,u_n(t)) dt} > 1/10 \min_{x>2^{n-1}}{f(t,x)}$
goes to infinity with $n$;
the $L^2$ norm of $v_n \in K$ is unbounded,
and we are done with this case.
Otherwise, let $t'_n$ be the largest value in $I_n$
for which $u_n(t'_n) = 2^{n-1}$:
$\int_{t'_n}^{t_n}{v_n(t) dt} = u_n(t_n) - u_n(t'_n) +
\int_{t'_n}^{t_n}{f(t,u_n(t)) dt} > 2^{n-1}$
for sufficiently large $n$ and again we have a contradiction.
\qed

{\nobf Remark: }
There are simple a priori estimates yielding properness
in the autonomous $H$ case.
Also, in the $C$ case (even for non-autonomous $f$)
easy estimates obtain Proposition 1.4 without invoking Theorem 1.2.
The analogous proof in the general $H$ case appears to be
considerably more elaborate
and we preferred making use of the more geometric Theorem 1.2.

The results above can be used to provide a simple definition
of topological degree for the operator $F$ in the case when $f$ is proper:
$$\deg F = \deg \Fbd =
\sum_{w \in \Fbd^{-1}(v)}{\sgn \phi_2(w)},$$
where $v$ is an arbitrary regular value of $F$.
As usual, the right hand side does not depend on the choice of $v$:
it is the degree of $\upsilon \mapsto \phi(\tilde v,\upsilon)$,
a proper function from $\RR$ to $\RR$.
From the behaviour of $\phi$ at infinity (proof of Proposition 1.4),
$$\deg F =
\sgn\left({\lim_{x \to +\infty}{f(t,x)}}\right) -
\sgn\left({\lim_{x \to -\infty}{f(t,x)}}\right).$$

Adapted coordinates give another simple characterization
of the critical set:
$(\tilde v, \overline u)$ is a critical point of
$\Fbd$ if and only if
$D_2\phi(\tilde v, \overline u) = 0$.
Equivalently, $u$ is a critical point of $F$ if and only if
$D\Phi(u) \wz = 0$ where $\wz$
is the tangent vector to $\tau_u$ at $u$ given by
the pull-back $\wz_u = (D\Psi(u))^{-1}({\bf 1}(\Psi(u)))$
(${\bf 1}$ is the vertical vector field
consisting of the constant function $1$ at each point).

{\nobf Lemma 1.5: }
{\sl
Given $u \in B^1$ and $m \in \RR$,
there is a unique $\alpha \in \RR$ such that the equation
$$\omega' + D_2f(t,u(t)) \omega = \alpha,\eqno{(\ast)}$$
has a (unique) periodic solution $\omega$ of average $m$.
The function $\wz$ is the only such $\omega$ of average 1;
furthermore, $\wz$ is strictly positive.
}

{\nobf Proof: }

The first claim follows from either solving $(\ast)$
or from arguments in the proof of Theorem 1.2.
For the element $u_a$ of average $a$ in $\tau_u$,
$$u_a'(t) + f(t, u_a(t)) - \int f(t, u_a(t)) dt = \tilde v.$$
Differentiating in $a$ and setting $\wz = {\partial \over \partial a}u_a$,
the equation $(\ast)$ for $\wz$ follows.
Since $\int u_a = a$, $\int \wz = 1$.
The graphs of $\wz$ and the constant function 0 do not cross,
implying positivity of $\wz$.
\qed

The following lemma introduces yet another characterization
of the critical set $S_1F$ and, under generic hypothesis,
establishes convenient transversality properties
of these characterizations.
This lemma will be essential for the more detailed study
of singularities of $F$ in the next section.

{\nobf Lemma 1.6: }
{\sl Let $\Sigma_a, \Sigma_b, \Sigma_c: B_1 \to \RR$ be given by
$$\eqalign{
\Sigma_a(u) &= \int D_2f(t, u(t)) dt,\cr
\Sigma_b(u) &= \int D_2f(t, u(t)) w_u(t) dt,\cr
\Sigma_c(u) &= \int D_2f(t, u(t)) \wz_u(t) dt.\cr}$$
Then the three $\Sigma$'s differ
by strictly positive smooth multiplicative factors.
In particular, $S_1F$ is the zero level of each of these functionals and
if 0 is a regular value of any of the functionals,
it is a regular value of all of them.}

{\nobf Remarks: }
\meti{1.}{For an open dense set of functions $f$,
0 is a regular value of $\Sigma_a$.
Indeed, taking derivatives as usual, 0 is a singular value if and only if
there is $u \in B^1$ with $\int D_2f(t, u(t)) dt = 0$
and $D_2D_2f(t, u(t)) = 0$ for all $t \in \SS^1$.}
\meti{2.}{In the autonomous case, 0 is a singular value of $\Sigma_a$
if and only if $D_2f$ has a double root.}

{\nobf Proof: }
Here, $P_i$ stands for a smooth strictly positive function.
From the expression for $w$ in Lemma 1.1,
$$w(t) = P_1(u) e^{-\int_0^t ( D_2f(s, u(s)) - \Sigma_a(u) ) ds},$$
Thus,
$$ \Sigma_b(u) = P_1(u) \int_0^1 D_2f(t, u(t))
e^{-\int_0^t ( D_2f(s, u(s)) - \Sigma_a(u) ) ds} dt. $$
On the other hand,
$$\eqalign{0 &= \int_0^1 {d \over dt} 
e^{-\int_0^t ( D_2f(s, u(s)) - \Sigma_a(u) ) ds} dt \cr
&= \int_0^1 (\Sigma_a(u) - D_2f(t, u(t)))
e^{-\int_0^t ( D_2f(s, u(s)) - \Sigma_a(u) ) ds} dt \cr}$$
whence
$$\eqalign{
\Sigma_b(u)
&= P_1(u) \Sigma_a(u) \int_0^1
e^{-\int_0^t ( D_2f(s, u(s)) - \Sigma_a(u) ) ds} dt \cr
&= P_2(u) \Sigma_a(u).\cr}$$

Integrate from 0 to 1 the differential equation describing $\wz$
to obtain $\alpha = \Sigma_c(u)$.
Solving the equation, we have
$$\wz(t) = \wz(0) e^{-\int_0^t D_2f(s, u(s)) ds}
+ \Sigma_c(u) e^{-\int_0^t D_2f(s, u(s)) ds}
\int_0^t e^{-\int_0^s D_2f(r, u(r)) dr} ds.$$
From $\wz(1) = \wz(0)$, we obtain
$$\wz(0) ( 1 - e^{-\Sigma_a(u)} ) =
\Sigma_c(u) e^{-\Sigma_a(u)}
\int_0^1 e^{-\int_0^s D_2f(r, u(r)) dr} ds.$$
Since $(1 - e^{-x})/x > 0$,
$$\Sigma_c(u) = P_3(u) \Sigma_a(u),$$
and we are done.
\qed

From now on, we shall always assume that $f$ is generic
in the sense that 0 is a regular value of $\Sigma_a$;
further generic properties will be required of $f$
in Section 2 where we study in detail the singularities of $F$.

We shall later want to use the simpler $w$ instead of $\wz$:
the following preparatory lemma allows for this interchange.

{\nobf Lemma 1.7: }
{\sl The vector fields $w$ and $\wz$ are positive multiples
of each other on $S_1F$. Furthermore, given $u \in S_1F$
there is a neighborhood $U_u \subseteq B^1$ of $u$ where we can write
$$w = a_1\wz + \Sigma_c z_1,\qquad
\wz = a_2 w + \Sigma_b z_2$$
for smooth real functions $a_i: U_u \to \RR$
and smooth vector fields $z_i$.}

{\nobf Proof:}

The first claim follows directly from the formulae for $w$ and $\wz$
when restricted to $S_1F$ (where $\Sigma_a = \Sigma_c = 0$).
The displayed equations are consequences of the regularity
of $\Sigma_b$ and $\Sigma_c$ at $S_1F$.
\qed

{\nobf Remark: }Actually, from results in [MST], $U_u$ in the statement
can be taken to be the whole space $B^1$.

It turns out that the global geometry of $S_1$
is very simple, as we shall see in Corollary 1.9.
We need some preparation to state the key ingredient, Theorem 1.8.

Let $M$ be a smooth compact manifold equipped
with a unit measure $\mu$.
Given a continuous function $g_k: M \times \RR \to \RR^k$,
define $G_k: B^1 \to \RR^k$
to be the average of the related Nemytski{\u\i} operator:
$G_k(v) = \int_M g_k(m, v(m)) d\mu$.
We request that $g_k$ admits continuous
partial derivatives of all orders
with respect to the second variable,
whence $G_k$ is smooth.

Let $\Pi_i: \RR^k \to \RR^i$ be the projection
to the first $i$ coordinates.
We say $0$ is a {\sl strong regular value} of $G_k$
if it is a regular value of the composition $G_i = \Pi_i \circ G_k$
for all $i$, $1 \le i \le k$.

{\nobf Theorem 1.8: {\rm [MST]}}
{\sl
Assume 0 to be a strong regular value of $G_k$.
Then the levels $Z_i$ ($1 \le i \le k$) are contractible manifolds.
Furthermore, there is a global homeomorphism $\Xi$
of $B^1$ taking each $Z_i$
to a closed linear subspace of codimension $i$;
$\Xi$ can be chosen to be a diffeomorphism if $B^1 = H^1$.
}

The contractibility of the levels $Z_i$ essentially implies
geometric triviality because of infinite dimension:
recall that two infinite dimensional separable Hilbert manifolds
are diffeomorphic if their homotopy groups coincide ([Ku])
and that all infinite dimensional separable Banach spaces 
are homeomorphic ([Ka]).

In the next corollary we have $k=1$;
Theorem 1.8 in its generality will be convenient in Section 3.

{\nobf Corollary 1.9: }
{\sl
Assume that 0 is a regular value of $\Sigma_a$.
Then $S_1$ is connected and contractible.
Furthermore,
there is a global homeomorphism $\Xi$ of $B^1$ taking $S_1$
to a closed linear subspace of $B^1$ of codimension $1$;
$\Xi$ can be chosen to be a diffeomorphism if $B^1 = H^1$.
}

\bigbreak

\bigbreak

{\noindent\bigbf 2. Morin theory}

\nobreak\smallskip\nobreak

Morin classified generic singularities of
functions from $\RR^n$ to $\RR^n$ whose derivative
has kernel of dimension 1 ([M]).
The first step in Morin's proof makes use of the implicit function
theorem to write such a singularity at the origin
in {\sl adapted coordinates}, i.e., in the form
$$(\xx, y) \mapsto (\xx, \mu(\xx, y)), \qquad
\xx = (x_1, \ldots, x_{n-1}) \in \RR^{n-1},\; y \in \RR,$$
after composing with suitable diffeomorphisms in the neighborhoods
of zero in both domain and image.
Morin's central result is that such singularities
are classified by their {\sl order}:
a {\sl Morin singularity} of order $k$
is a point $(\xx,y)$ for which
\item{(a)}$D_2\mu(\xx,y) = \cdots = D_2^k\mu(\xx,y) = 0$,
\item{(b)}$D_2^{k+1}\mu(\xx,y) \ne 0$,
\item{(c)}the Jacobian $D(D_2\mu, \ldots, D_2^{k-1}\mu)(\xx,y)$ is surjective.

Set $S_k = \{(\xx,y) | D_2\mu(\xx,y) = \cdots = D^k_2\mu(\xx,y) = 0\}$.
Thus, $S_0$ is the domain,
$S_1$ is the critical set $\{(\xx,y) | D_2\mu(\xx,y) = 0\}$
(consistently with previous notation).
Also, a Morin singularity of order $i$ belongs to $S_k$ if and only if
$i \ge k$.
In a neighborhood of a Morin singularity,
the sets $S_k$ stratify the domain:
the sets are nested and $S_i$ is a submanifold of codimension $i$.
Notice that a point $(\xx,y) \in S_k - S_{k+1}$
is a Morin singularity (of order $k$) only if condition (c) above holds.

Composing by appropriate diffeomorphisms in the domain and image,
a Morin singularity in dimension $n$ and order $k$
acquires the normal form
$$(x_1, \ldots, x_{n-1}, y) \mapsto
(x_1, \ldots, x_{n-1}, y^{k+1} + x_1 y^{k-1} + \cdots + x_{k-1} y).$$
Morin singularities of order 1, 2, 3 and 4 are called, respectively,
folds, cusps, swallowtails and butterflies.

We shall need an equivalent classification
for singularities of functions between infinite-dimensional spaces.
Let $G: Z_1 \to Z_2$ be a smooth map between Banach spaces
so that $DG(z_0)$ is Fredholm operator of index 0 and
kernel of dimension 1.
Again, after changes of variables in the domain and image,
we may assume $G$ near $z_0$ to be written in adapted coordinates as
$$\eqalign{
G: X \times \RR &\to X \times \RR,\cr
(\xx, y) &\mapsto (\xx, \mu(\xx, y))\cr}$$
and if conditions (a), (b) and (c) above hold,
the same normal form applies for an appropriate splitting
$X = \RR^{k-1} \oplus X'$---%
we then call $z_0$ a Morin singularity of order $k$.
The proof of this last fact follows Morin's ([M]),
making use of a parameterized version of
Malgrange's preparation theorem,
the parameter taking values in a Banach space (see [CDT]).

We already saw in the previous section
that the composition $\Fbd = F \circ \Psi^{-1}$ is in adapted coordinates:
$$\eqalign{
\Fbd: \tilde{B^0} \oplus \langle 1\rangle
&\to \tilde{B^0} \oplus \langle 1\rangle.\cr
(\tilde v, \overline u) &\mapsto
(\tilde v, \overline v = \phi(\tilde v, \overline u))}$$
From the previous paragraph,
a point is a Morin singularity of $\Fbd$ (or $F$) of order $k$ if and only if
conditions (a), (b) and (c) hold for $\mu$ replaced by $\phi$.
This criterion, however, can not be used directly since we have no formula
for $\phi$; we rephrase it in terms of $\Phi = \phi \circ \Psi$ and $w$.
Following the usual notation, we write
$w\xi$ for the Lie derivative $D\xi(u) \cdot w_u$.
The following result shows that we may substitute $w$ for $\wz$
(alternative generators for the kernel of $DF$ over $S_1$)---%
a fact which in finite dimension would be unsurprising.

{\nobf Proposition 2.1: }
{\sl The point $u \in B^1$ is a Morin singularity of order $k$
for $F$ if and only if:
\item{\rm (a)}$w\Phi(u) = \cdots = w^k\Phi(u) = 0$,
\item{\rm (b)}$w^{k+1}\Phi(u) \ne 0$,
\item{\rm (c)}$D(w\Phi, \ldots, w^{k-1}\Phi)(u)$ is surjective.}

{\nobf Proof: }
Consider in $\tilde{B^0} \oplus \langle 1\rangle$
the constant vertical vector field ${\bf 1}$,
consisting of the constant function $1$ at each point.
In the notation we just introduced,
$D_2\xi = {\bf 1}\xi$.
In terms of the pull-back $\wz(u) = (D\Psi(u))^{-1}({\bf 1}(\Psi(u)))$,
the conditions for $u$ to be a Morin singularity of order $k$ are:
\item{(a')} $\wz\Phi(u) = \cdots = \wz^k\Phi(u) = 0$,
\item{(b')} $\wz^{k+1}\Phi(u) \ne 0$ 
\item{(c')} $D(\wz\Phi, \ldots, \wz^{k-1}\Phi)(u)$ is surjective.

We are left with showing that we can substitute $\wz$ by $w$
in these conditions.
Notice first that $w\Phi = \Sigma_b$ and $\wz\Phi = \Sigma_c$,
proving the case $k=0$ (regular points).
From now on, we assume $u \in S_1F$ and write,
making use of Lemma 1.7, $w = a_1 \wz + (\wz\Phi) z_1$
and $\wz = a_2 w + (w\Phi) z_2$
in a small neighborhood $U_u$ of $u$.

For each $k$, the ideals in $C^\infty(U_u,\RR)$
generated by $w\Phi, \ldots w^k\Phi$ and
$\wz\Phi, \ldots \wz^k\Phi$ are equal.
Assuming by induction that the result holds for $k-1$,
$$\eqalign{
w^k\Phi &= w(w^{k-1}\Phi)\cr
&= (a_1 \wz + (\wz\Phi) z_1)(b_1 \wz\Phi + \cdots + b_{k-1} \wz^{k-1}\Phi)\cr
&= a_1 (\wz b_1)(\wz\Phi) + a_1 b_1 \wz^2\Phi + \cdots 
a_1 (\wz b_{k-1})(\wz^{k-1}\Phi) + a_1 b_{k-1} \wz^k\Phi +  \cr
&\qquad\qquad(\wz\Phi)z_1(b_1 \wz\Phi + b_{k-1} \wz^{k-1}\Phi)}$$
which is clearly in the ideal with generators
$\wz\Phi, \ldots \wz^k\Phi$, proving one inclusion;
the opposite inclusion is analogous.

The equality of the two ideals with $k$ generators implies
$$\eqalign{
w^k\Phi &= b_k \wz^k\Phi + \cdots + b_1 \wz\Phi,\cr
\wz^k\Phi &= c_k w^k\Phi + \cdots + c_1 w\Phi,\cr
}$$
for smooth functions $b_i$ and $c_i$
where $b_k$ and $c_k$ are non-zero.
The equivalence between the conditions (a) and (a') or (b) and (b') is clear.
The third equivalence follows from repeated use of the simple fact
that the spans of $D(g_1(u) g_2(u))$ and $D(g_1(u) g_2(u) + \alpha(u)g_1(u))$
coincide for points $u$ such that $g_1(u) = g_2(u) = 0$
($\alpha$ being a smooth real function).
\qed

{\nobf Proposition 2.2: }
{\sl For
$$\eqalign{
\Sigma_1(u) &= \ints D_2f(t,u(t)) dt,\cr
\Sigma_2(u) &= \ints D^2_2f(t,u(t)) w(t) dt,\cr
\Sigma_3(u) &= \ints D^3_2f(t,u(t)) w^2(t) dt,\cr
\Sigma_4(u) &= \ints D^4_2f(t,u(t)) w^3(t) - 
2 D^3_2f(t,u(t)) w^2(t) \left(\int_0^t D^2_2f(s,u(s)) w(s) ds\right) dt,\cr
\Sigma_5(u) &= \ints D^5_2f(t,u(t)) w^4(t) -
5 D^4_2f(t,u(t)) w^3(t) \left(\int_0^t D^2_2f(s,u(s)) w(s) ds\right) + \cr
&\qquad\qquad
5 D^3_2f(t,u(t)) w^2(t) \left({\int_0^t D^2_2f(s,u(s)) w(s)}\right)^2 dt,\cr}$$
we have, for $k = 1,\ldots, 5$,
$$S_k = \{ u \in B^1 | \Sigma_i(u) = 0, i = 1,\ldots,k \}.$$
Furthermore, for $k = 1, \ldots, 4$,
$u$ is a Morin singularity of order $k$ if and only if
\item{\rm (a)}$\Sigma_i(u) = 0$, $i = 1,\ldots,k$,
\item{\rm (b)}$\Sigma_{k+1}(u) \ne 0$,
\item{\rm (c)}the derivative $D\Sigma(u)$ of the function 
$$\eqalign{ \Sigma: B^1 &\to \RR^{k-1}\cr
u &\mapsto (\Sigma_1(u),\ldots,\Sigma_{k-1}(u))\cr}$$ is surjective.
}

{\nobf Proof: }
From Proposition 2.1, we want to compute $w^k\Phi(u)$.
The expressions for $\Sigma_k(u)$ follow from repeated integration
by parts, discarding elements in the ideal generated by
$w\Phi(u),\ldots,w^{k-1}\Phi(u)$ and non-zero multiplicative factors.
\qed

In particular, $u$ is a fold point if and only if
$\Sigma_1(u) = 0$ and $\Sigma_2(u) \ne 0$.
Also, $u$ is a cusp point if and only if
$\Sigma_1(u) = \Sigma_2(u) = 0$, $\Sigma_3(u) \ne 0$
and $D\Sigma_1(u) \ne 0$.
Clearly, $D\Sigma_1(u) \cdot v = \ints D_2^2f(t,u(t)) v(t) dt$,
and hence $D\Sigma_1(u) = 0$ if and only if
the function $D_2^2f(t,u(t))$ is identically 0.

% Funcoes batida

There is a simple relationship between Morin singularities
and the {\sl return map}.
Given $f: \SS^1 \times \RR \to \RR$ and $v \in B^0$,
the return map $\rho_v: I \to \RR$, $I \subseteq \RR$,
sends $u(0)$ to $u(1)$ if $u: [0,1] \to \RR$ satisfies
$$u'(t) + f(t,u(t)) = v(t).$$
Here $I$ is the maximal domain, i.e.,
the set of initial conditions such that the solution $u$
extends to $t = 1$.

{\nobf Proposition 2.3: }
{\sl For $u \in B^1$ and $v = F(u)$, $\rho_v(u(0)) = u(0)$.
Also, $u \in S_1$ if and only if $\rho_v'(u(0)) = 1$
and $u \in S_k$ if and only if $\rho_v'(u(0)) = 1$
and $\rho_v^{(i)}(u(0)) = 0$ for $i = 2, \ldots, k$.}

{\nobf Proof: }
Consider the fibre $\tau_u$ through $u$.
Let $u_a$ be the element of average $a$ of $\tau_u$ and $u = u_{a_0}$:
by Lemma 1.3, $\tau_u$ is smoothly parametrized by $a$.
Let $g(a) = \Phi(u_a)$, the average of $F(u_a)$;
clearly, $g^{(k)}(a) = \wz^k\Phi(u_a)$ and therefore
$u \in S^k$ if and only if $g^{(i)}(a_0) = 0$ for $i = 1, \ldots, k$.
Let $u^c$ be the element of $\tau_u$ with $u^c(0) = c$:
by Lemma 1.3, $\tau_u$ is also parametrized by $c$.
Let $h(c) = \Phi(u^c) - \overline v = F(u_c) - F(u)$: by the chain rule,
$u \in S^k$ if and only if $h^{(i)}(u(0)) = 0$ for $i = 1, \ldots, k$.

For $c, b \in \RR$, let $\beta(c,b) = \rho_{v+b}(c)$:
that is, if $U(c,b,t)$ satisfies
$$D_3U(c,b,t) + f(t, U(c,b,t)) = v(t) + b,\qquad U(c,b,0) = c\eqno{(\ast)}$$
we have $\beta(c,b) = U(c,b,1)$.
The periodicity of $u^c$ yields $\beta(c, h(c)) = c$.
Points in the curve $(c, h(c))$ thus correspond to points in $\tau_u$
and the largest $k$ for which $u \in S_k$ is the order of contact
between this curve and the horizontal axis $(c, 0)$
at the common point $(u(0), 0)$. 

Differentiating $(\ast)$,
$$D_3D_2U(u(0),0,t) + D_2f(t,u(t)) D_2U(u(0),0,t) = 1,
\qquad D_2U(u(0),0,0) = 0$$
and, by explicitly solving for $D_2U$, we obtain
$D_2U(u(0),0,1) = D_2\beta(u(0), 0) > 0$.
Thus, $G(c,b) = (c, \beta(c,b))$
is a local diffeomorphism near $(u(0), 0)$
taking the curve $(c,h(c))$ to the diagonal $(c,c)$
and the horizontal axis to $(c,\rho_v(c))$.
The order of contact between the curves
is preserved by $G$ and we are done.
\qed

\bigbreak

\bigbreak

{\noindent\bigbf 3. The autonomous case}

\nobreak\smallskip\nobreak

This section is dedicated to a number of special properties
of the {\sl autonomous case}, when $f$ depends on $x$ only.
Notice that this implies that $f$ is tame.

The sets $S_k$ are described by the rather complicated formulae $\Sigma_k$;
in the autonomous case it is convenient to consider
the {\sl simplified operator}
$$\eqalign{
\hat F:B^1 &\to B^0\cr
u &\mapsto u' + \int f(u(t)) dt\cr}$$
whose critical strata $\hat S_k$
are far easier to handle but still convey significant information
about $F$ and $S_k$.

Since $\hat F$ is already given in adapted coordinates,
straightforward application of Morin's characterization
obtains
$$\hat S_k = \{v \in B^1 | \ints\hat\gamma_k(v(t)) = 0\},$$
where $\hat\gamma_k(x)$ is the $k$-dimensional vector
$(f'(x),\ldots,f^{(k)}(x))$
Notice that, from Lemma 1.6, $\hat S_1 = S_1$.
However, $\hat S_2$ is usually different from $S_2$:
the following lemma relates both sets and is a key ingredient in Section 5.

{\noindent\bf Lemma 3.1:}
{\sl
Let $f: \RR \to \RR$ be a smooth function,
$F: B^1 \to B^0$ be the operator $(F(u))(t) = u'(t) + f(u(t))$.
Then there exists a global diffeomorphism of $B^1$
taking $S_1$ to itself and $S_2$ to $\hat S_2$.
}

{\noindent\bf Proof:}
The diffeomorphism has the form
$u \mapsto v = u \circ \alpha$,
with inverse $v \mapsto u = v \circ \beta$, 
where $\alpha$ and $\beta$ are orientation preserving
$C^1$ diffeomorphisms of $\SS^1$ fixing 0:
such compositions take functions in $B^1$ to functions in $B^1$.
Clearly, $\beta = \alpha^{-1}$.
We need the following characterizations,
which follow easily from Lemma 1.6 and Proposition 2.2:
$$\eqalign{
S_1 &=
\left\{{u \in B^1 \;\Big|\; \ints f'(u(t)) w(t) dt = 0}\right\},\cr
S_2 &= \left\{{u \in B^1 \;\Big|\; \ints f'(u(t)) w(t) dt =
\ints f''(u(t)) w(t) dt = 0}\right\}.\cr}
$$

We first obtain $\beta$ from $u$:
$$\beta(s) = {{\int_0^s w(\sigma) d\sigma} \over
{\int_0^1 w(\sigma) d\sigma}}$$
is clearly a $C^1$ diffeomorphism ($w > 0$).
For any continuous function $g$,
a change of variables gives
$${\int_0^1 g(u(s)) w(s) ds = 0} \quad \Longleftrightarrow
\quad {\int_0^1 g(v(t)) dt = 0},\eqno{(\ast)}$$
where $u = v \circ \beta$.
Thus, if $u \in S$ then $v \in S$ and
if $u \in C$ then $v \in \hat S_2$.
The smooth dependence of $\beta$ on $u$ is obvious.

To show invertibility of the map $u \mapsto v$,
we obtain $\alpha$ from $v$.
If $\alpha$ satisfies
$$\alpha'(t) = {1 \over {A - 1 \int_0^tf'(v(\tau))d\tau}},
\quad \alpha(0) = 0$$
for an arbitrary positive constant $A$,
standard algebra shows that
$$\alpha'(t) = {1\over A} \exp\left({
\int_0^t f'(v(\tau)) \alpha'(\tau) d\tau}\right),\quad
\alpha(0) = 0.\eqno{(\dag)}$$
This again implies the equivalence $(\ast)$
(where, of course, $v = u \circ \alpha$) and hence
$v \in S$ (resp., $\hat S_2$) implies $u \in S$ (resp., $C$)
provided $\alpha(1) = 1$.
We have to show that for each $v$ there is a unique $A$
with $\alpha(1) = 1$ and that the dependence of $A$ on $v$ is smooth.

Let $h(t) = \int_0^t f'(v(\tau)) d\tau$;
$h$ is $C^1$ with $h(0) = 0$.
The function $\alpha$ is defined in $[0,1]$ if $A > \max_t h(t)$.
From $$\alpha(1) = \int_0^1 { dt \over A - h(t) },$$
the derivative of $\alpha(1)$ with respect to $A$ is strictly negative.
When $A$ tends to infinity, $\alpha(1)$ becomes small
and when $A$ approaches $\max_t h(t)$, $\alpha(1)$ tends to infinity.
This settles existence and uniqueness of the required $A$;
smoothness follows from the implicit function theorem
applied to the smooth function $(v,A) \mapsto \alpha(1)$
and the fact that the derivative with respect to $A$ is not zero.

It remains only to show that the two smooth maps constructed above
are the inverse of each other.
Consider the sequence of maps
$$v {\buildrel \alpha \over \longrightarrow} u = v \circ \alpha^{-1}
{\buildrel \beta \over \longrightarrow} \tilde v = u \circ \beta^{-1}.$$
We have, for positive constants $C_1$ and $C_2$,
$$\eqalignno{
\beta'(s) &= C_1 w(s) \cr &=
C_1 \exp\left({-\int_0^s f'(v(\alpha^{-1}(\sigma))) d\sigma }\right) \cr
&{\buildrel{\scriptstyle\tau = \alpha^{-1}(\sigma)}\over{=}}
C_1 \exp\left({-\int_0^{\alpha^{-1}(s)}
f'(v(\tau)) \alpha'(\tau) d\tau}\right)
\hfill{(\hbox{from }\tau = \alpha^{-1}(\sigma))}\cr
&{\buildrel{(\dag)}\over{=}} {C_2 \over \alpha'(\alpha^{-1}(s))} \cr
&= C_2 (\alpha^{-1})'(s) \cr
}$$
and, since $\alpha(0) = \beta(0) = 0$ and $\alpha(1) = \beta(1) = 1$,
it follows that $\alpha^{-1} = \beta$ and $\tilde v = v$.
Similarly, for
$$u {\buildrel \beta \over \longrightarrow} v = u \circ \beta^{-1}
{\buildrel \alpha \over \longrightarrow}
\tilde u = v \circ \alpha^{-1},$$
we have
$$\eqalignno{
\alpha'(t) &=
C_3 \exp\left({\int_0^t f'(v(\tau)) \alpha'(\tau) d\tau}\right) \cr
&{\buildrel{\scriptstyle\tau = \alpha^{-1}(\sigma)}\over{=}}
C_3 \exp\left({\int_0^{\alpha(t)}
f'(v(\alpha^{-1}(\sigma))) d\sigma }\right) \cr
&= {C_3 \over w(\alpha(t))}\cr
&= {C_4 \over \beta'(\alpha(t))},\cr
}$$
hence $(\beta \circ \alpha)'$ is a constant
and the result follows.
\qed

As with $S_1 = \hat S_1$,
the global geometry of the sets $\hat S_k$ is very simple,
as we learn from the following application of Theorem 1.8.
We say that the simplified operator $\hat F$ is {\sl $k$-regular}
if 0 is a strong regular value of
$$\eqalign{
\hat\Sigma: B^1 &\to \RR^k.\cr
u &\mapsto (\hat\Sigma_1,\ldots,\hat\Sigma_k)\cr}$$

{\nobf Corollary 3.2: }
{\sl
Assume $\hat F$ is $k$-regular.
Then there is a global homeomorphism $\Xi$ of $B^1$ taking each $\hat S_i$
($1 \le i \le k$)
to a closed linear subspace of $B^1$ of codimension $i$;
$\Xi$ can be chosen to be a diffeomorphism if $B^1 = H^1$.
}

Combining Lemma 3.1 and Corollary 3.2 we have:

{\nobf Corollary 3.3: }
{\sl
Assume $\hat F$ is 2-regular.
Then there is a global homeomorphism $\Xi$ of $B^1$ taking each $S_i$,
$i = 1, 2$,
to a closed linear subspace of $B^1$ of codimension $i$;
$\Xi$ can be chosen to be a diffeomorphism if $B^1 = H^1$.
}

\bigbreak

It seems hard to give an operational criterion
to decide for larger $k$ even whether $S_k$ is non-empty.
We now present a partial criterion.

{\noindent\bf Definition 3.4:}
{\sl 
Let $f:\RR \to \RR$ be a smooth function;
$f$ is said to be {\rm $k$-good} if
$\hat\gamma_k$ never vanishes and 
the image of any open interval by $\hat\gamma_k$
is not contained in a hyperplane through the origin in $\RR^k$.
}

Generic smooth functions are $k$-good,
as well as generic polynomials of fixed degree at least $k$.
It is easy to see that if $f$ is $k+1$-good,
then the simplified operator $\hat F$ is $k$-regular.

{\noindent\bf Lemma 3.5:}
{\sl
Let $f:\RR \to \RR$ be a $k$-good function
and $F$ be the related operator.
Then $\hat S_k \ne \emptyset$ if and only if $0 \in \RR^k$
is in the interior of the convex hull of the image of $\hat\gamma_k$.
Also, $\hat S_k \ne \emptyset$ implies $S_k \ne \emptyset$.
}

In section 6, we give an example of a polynomial of degree 4
(and hence for which $\hat S_4 = \emptyset$)
having $S_4 \ne \emptyset$.

\goodbreak

{\noindent\bf Proof:}

Clearly, if $\hat S_k \ne \emptyset$,
$0$ is in the convex hull of the curve $\hat\gamma_k$
(but not on the curve itself).
If 0 is on the boundary of the convex hull,
by a standard support theorem ([G], pg. 12),
the image of $\hat\gamma_k$ is contained in a closed half-space
defined by $\nu(p) \ge 0$ for some linear functional $\nu$.
Since $f$ is $k$-good,
for any non-constant $v \in \hu$,
$\nu(\ints\hat\gamma_k(v)) > 0$ and $\hat S_k = \emptyset$.

Conversely, assume 0 in the interior of
the convex hull of the image of $\hat\gamma_k$.
By Steinitz's theorem ([G]),
there are points $\hat\gamma_k(x_j)$, $j = 0,\ldots,2k-1$
such that 0 is in the interior of their convex hull.
For $\epsilon \in (0,1)$ and
$a_j \ge 0$, $j = 0,\ldots,2k-1$, $\sum a_j = 1 - \epsilon$,
consider a smooth function $\upsilon_{\epsilon,a}$ of period 1,
defined as follows.
We split the domain $[0,1]$ into intervals
$I_0,J_0,\ldots,I_{2k-1},J_{2k-1}$
of lengths $a_0,\epsilon/{2k},\ldots,a_{2k-1},\epsilon/{2k}$;
inside $I_j$, $\upsilon_{\epsilon,a}$ is constant equal to $x_j$
and inside each $J_j$, $\upsilon_{\epsilon,a}$ is the appropriate
affine transformation of a fixed smooth arc joining two steps.
As $\epsilon$ tends to 0,
$\upsilon_{\epsilon,a}$ approaches a step function.
Let $\phi(\epsilon,a) = \ints\hat\gamma_k(\upsilon_{\epsilon,a})$:
this function is affine in $\epsilon$ and $a$
(i.e., linear plus constant)
and to show that $\hat S_k$ is non-empty,
we need to find a zero of $\phi$.
The function $\phi$ extends continuously to $\epsilon=0$
and $0$ is then in the interior of the image of the simplex
spanned by the $a_j$'s:
there exists therefore a straight segment parametrized
by small positive values of $\epsilon$ on which $\phi$ is zero
and $\hat S_k$ is thus non-empty.
More, for a fixed small $\epsilon_0$,
take a $k$-subspace $V$ of $\RR^{2k}$
such that, for $a \in V$, $\phi(\epsilon_0,a)$ is surjective.
The image under $\phi$ of a small sphere around the origin in $V$
is some ellipsoid containing the origin in its interior.

Now, assume $\hat S_k \ne \emptyset$.
Use the space $V$ and the function $\phi$
to obtain $r > 0$ and a function $H: \BB^k \to B^1$ with
$$ \int \hat\gamma_k(H(s)(t)) dt = r s $$
for $s \in \BB^k$ where $\BB^k \subset \RR^k$ is the unit ball.
Define the {\it $N$-replicator} to be the isomorphism
$R_N: B^i \to R_N(B^i) \subset B^i$, $(R_N(u))(t) = u(Nt)$, $i = 0, 1$. 
Clearly, $(R_N(u))' = N R_N(u')$.
We claim that given $\epsilon > 0$ there exists $N$
such that
$$ \left| \left( w\Phi(R_N(H(s))),\ldots, w^k\Phi(R_N(H(s))) \right)
- rs \right| < \epsilon $$
for all $s \in \BB^k$ and the proof is completed
by a standard degree theory argument.

At this point it is convenient to make explicit the dependence
of $w = w(u)$ in terms of $u$.
From Proposition 1.1,
% $$ (w(u))'(t) + f'(u(t)) w(u)(t) = \lambda w(u)(t), \quad
$$ (w(R_N(u)))'(t) + f'(R_N(u(t))) w(R_N(u))(t) =
\lambda w(R_N(u))(t), \quad
\lambda = \int f'(R_N(u)(t)) dt $$
and therefore $\lambda$ is the same for $u$ and $R_N(u)$.
Define $w_N(u)$ by $R_N(w_N(u)) = w(R_N(u))$ so that
$$ (w_N(u))'(t) + {{f'(u(t))}\over{N}} w_N(u)(t) =
{{\lambda}\over {N}} w_N(u)(t) $$
and, from the formula for $w$ in Proposition 1.1,
$w_N(u)(t) = (w(u)(t))^{(1/N)}$.
Since $\Phi(u) = \int f(u(t)) dt$ we have
$\Phi(R_N(u)) = \Phi(u)$ and
$$ \left( w\Phi(R_N(u)),\ldots, w^k\Phi(R_N(u)) \right) 
= \left( w_N\Phi(u),\ldots, w_N^k\Phi(u) \right). $$
The sequence $(w_N)$ of vector fields tends to the constant vector field
$\fancyone$ (i.e., the constant function $1$ at every point $u$)
in the $C^n$-metric (for any $n$).
Also, 
$$ \left( \fancyone\Phi(R_N(u)),\ldots, \fancyone^k\Phi(R_N(u)) \right)  =
\int\hat\gamma_k(u(t)) dt, $$
proving the claim.
\qed

{\noindent\bf Remarks:}
For $k = 2$, $k$-goodness may be weakened to
$f'$ and $f''$ having no common zeros.

It also follows from the above argument that 
$D((w\Phi, \ldots, w^k\Phi) \circ R_N \circ H)$ tends to the identity matrix
when $N$ goes to infinity, establishing condition (c) in Proposition 2.1.
Condition (b) follows from the additional hypothesis of $f$ being $(k+1)$-good,
thus proving the existence of Morin singularities of order $k$.

\bigbreak

\bigbreak

{\noindent\bigbf 4. Some examples}

\nobreak\smallskip\nobreak

In this section, we describe the global geometry of $F: B^1 \to B^0$
for several special classes of smooth functions $f: \SS^1 \times \RR \to \RR$.
The first rather simple proposition illustrates the use of fibres
and adapted coordinates in three technically different scenarios.

{\nobf Proposition 4.1: }
{\sl
\item{\rm (a)}
If $f$ is proper and $D_2f(t,x) > 0$ then $F$ is a diffeomorphism.
\item{\rm (b)}
Assume $D_2f(t,x) > 0$. Let $f_\pm(t) = \lim_{x \to \pm\infty} f(t,x)$.
Then $F$ is a diffeomorphism from $B^1$ to the horizontal strip
$${\cal S} =
\{v \in B^0 | \ints f_-(t) dt < \ints v(t) dt < \ints f_+(t) dt\}.$$
\item{\rm (c)}
If $f$ is proper and strictly increasing in the second variable
then $F$ is a homeomorphism.}

{\nobf Proof: }
Notice that the hypotesis on each item imply that the nonlinearity $f$ is tame.

\meti{(a)}Recall that, by Theorem 1.2, $B^1$ is foliated by fibres,
$B^0$ is foliated by vertical lines
and there is a diffeomorphism between the space of fibres
and the space of vertical lines.
From the characterization of $S_1$ in Proposition 1.1,
we see that $F$ has no critical points.
Thus, $F$ takes each fibre strictly monotonically
to its related vertical line in $B^0$---%
by Proposition 1.4, $F$ is actually a diffeomorphism
from fibre to vertical line and the result follows.

\meti{(b)}
As in the previous item, fibres are bijectively taken
to open subintervals of vertical lines.
More explicitly, a function $u$ in the fibre
$(\Pi_{\tilde B^0} \circ F)^{-1}(\tilde v)$
is taken to $u' + f(t,u(t)) = \tilde v + \ints f(t,u(t)) dt$
and the extremes of the image of this fibre are
$\lim_{a \to +\infty} \ints f(t, u_a(t)) dt$ and
$\lim_{a \to -\infty} \ints f(t, u_a(t)) dt$, where
$u_a$ is the element of average $a$ in the fibre.
From Proposition 1.3,
$$\eqalign{
\lim_{a \to +\infty} \min_t u_a(t) &= +\infty,\cr
\lim_{a \to -\infty} \max_t u_a(t) &= -\infty.\cr}
$$
Thus,
$$\eqalign{
\lim_{a \to +\infty} \ints f(t, u_a(t)) dt &= \ints f_+(t) dt,\cr
\lim_{a \to -\infty} \ints f(t, u_a(t)) dt &= \ints f_-(t) dt,\cr}
$$
and the result follows.

\meti{(c)}By properness (Proposition 1.4),
it suffices to prove that $F$ is strictly increasing on each fibre
(notice that $F$ restricted to fibres may have critical points
where $\wz F = 0$).
Let $u_0 < u_1$ be elements of the fibre
$(\Pi_{\tilde B^0} \circ F)^{-1}(\tilde v)$ so that
$$F(u_i)(t) = u_i'(t) + f(t, u_i(t)) = \tilde v + \overline v_i.$$
Integrating in $t \in \SS^1$ and using the monotonicity of $f$,
we obtain $\overline v_0 < \overline v_1$, concluding the proof.
\qed

{\nobf Remarks: }
\meti{1.}
A more standard proof of (a), without making use of Theorem 1.2
(and the consequent fibre-sheet-adapted coordinates vocabulary),
could be as follows.
As before, $S_1 = \emptyset$ and
from Proposition 1.4, $F$ is proper.
Since $B^0$ is simply connected, 
by covering space theory, $F$ is a diffeomorphism.
Notice that this argument does not extend easily to the other items.
\meti{2.} Similar results and proofs hold if instead $D_2f < 0$.
\meti{3.} In item (b), if both $\ints f_+(t) dt$ and $\ints f_-(t) dt$ diverge,
then $F$ is a global diffeomorphism.
In particular, $F$ may be proper even if $f$ is not.

{\nobf Theorem 4.2: }
{\sl If $f$ is proper and $D_2^2f(t,x) > 0$ then $F$ is a global fold.}

We call an operator $G: B^1 \to B^0$ a {\sl global fold}
if there exist diffeomorphisms
$\Xi_1: B^1 \to \RR \times \tilde B^0$ and 
$\Xi_0: B^0 \to \RR \times \tilde B^0$
such that $(\Xi_0 \circ G \circ \Xi_1^{-1})(x, \tilde v) = (x^2, \tilde v)$,
for all $(x, \tilde v) \in \RR \times \tilde B^0$.
Similarly, we call $G$ a {\sl topological} global fold
if there exist homeomorphisms $\Xi_i$ as above.

{\nobf Proof: }
From $D_2^2f(t,x) > 0$, we conclude that $D_2(t,x)$ is strictly increasing
in $x$ for any fixed $t$ and hence that $\Sigma_a(u) = \ints D_2(t,u(t)) dt$
is strictly increasing on fibres.
Thus, each fibre contains a unique critical point $u_0$
and, for arbitrary $u_-$ and $u_+$ in the same fibre as $u_0$
satisfying $u_- < u_0 < u_+$, we have
$\Sigma_a(u_-) < 0$ and $\Sigma_a(u_+) > 0$
and, from Lemma 1.6,
$\wz\Phi(u_-) = \Sigma_c(u_-) < 0$ and $\wz\Phi(u_+) = \Sigma_c(u_+) > 0$
and the restriction of $\Phi$ to a fibre is a global fold from $\RR$ to $\RR$.
Thus, on each fibre, we have diffeomorphisms $\Xi_i$ as above
and the problem is whether such diffeomorphisms can be chosen
so as to depend smoothly on the fibre.

In adapted coordinates,
we must define $\xi_1: B^0 \times \RR \to \RR$
and $\xi_0: B^0 \times \RR \to \RR$
so that the vertical columns of the diagram below are diffeomorphisms
and the diagram commutes, i.e.,
$$\xi_0(\tilde v, \phi(\tilde v, \overline u)) =
(\xi_1(\tilde v, \overline u))^2.$$
We construct the $\xi_i$ explicitly.
For each $\tilde v$, let $a_{\tilde v}$ be the unique critical point
of $\overline u \mapsto \phi(\tilde v, \overline u)$.
Clearly, $a_{\tilde v}$ and its image
$b_{\tilde v} =  \phi(\tilde v, a_{\tilde v})$
depend smoothly on $\tilde v$.
Also, write $\phi(\tilde v, \overline u) - b_{\tilde v} =
(\overline u - a_{\tilde v})^2 g(\tilde v, \overline u)$.
From the previous paragraph, $g$ is a smooth positive function.
Set $\xi_0(\tilde v, \overline v) = \overline v - b_{\tilde v}$
and $\xi_1(\tilde v, \overline u) =
(\overline u - a_{\tilde v}) \sqrt{g(\tilde v, \overline u)}$.
\qed

{\nobf Remarks: }
\meti{1.}
McKean and Scovel ([McKS]) studied this scenario
with a different set of fibres for $B^1$ and its image under $F$.
Our choice of fibering $B^0$ by vertical lines
and inverting them under $F$ to get a fibration of $B^1$
is more helpful in our examples.
\meti{2.} Slight variations (as in Proposition 4.1) are possible
and can be handled similarly; we omit the tedious details.

In the autonomous case, theorem 4.2 admits a partial converse:

{\nobf Theorem 4.3: }
{\sl
Let $f: \RR \to \RR$ be a smooth function which is both 2- and 3-good,
with $\lim_{x \to \pm \infty} f(x) = +\infty$.
If 0 is in the interior of the convex hull of the image
of $\hat\gamma_2$ then $F$ has cusps.
Otherwise, $F$ is a (differentiable) global fold.}

Recall that $\hat\gamma_2(t) = (f'(t),f''(t))$.
Notice that once $F$ is known to have a cusp,
from the normal form we have image points
with three regular pre-images near the cusp;
since $\deg(F) = 0$, such points have at least one additional pre-image.

{\nobf Proof: }
If 0 is in the interior of the convex hull
of the image of $\hat\gamma_2$,
$S_2 \ne \emptyset$ by Lemma 3.5.
We must prove that some points $u$ in $S_2$ are differentiable cusps, i.e.,
satisfy $\Sigma_3(u) \ne 0$ and $D\Sigma_1(u) \ne 0$.
As remarked after Proposition 2.2, $D\Sigma_1(u) = 0$ only when $f''(u(t)) = 0$
for all $t$, which implies, given 2-goodness, that $u$ is constant
equal to a root of $f''$. Again from 2-goodness,
$f'$ and $f''$ have no common roots and
$\Sigma_1(u) = \ints f'(u(t)) dt \ne 0$,
thus $u$ is not in the critical set.
It remains to verify that we may choose $u$ so that $\Sigma_3(u) \ne 0$.
From 3-goodness, the curve $\hat\gamma_3(t) = (f'(t),f''(t),f'''(t))$
does not intersect the origin.
The convex hull of the image of $\hat\gamma_3$ meets the vertical axis
and must contain points $(0,0,A)$ distinct from the origin,
otherwise the image of $\hat\gamma_3$ would have to be contained
in a hyperplane, contradicting 3-goodness.
Imitating the proof of Lemma 3.5,
we may construct $u$ with
$\Sigma_1(u) = \Sigma_2(u) = 0$, $\Sigma_3(u) \approx A$,
which is the required cusp.

Conversely, assume 0 not to be in the interior of the convex hull
of the image of $\hat\gamma_2$.
Again from Lemma 3.5, all critical points are folds.
Notice that if $m$ is the minimum of $f$,
$f'(m) = 0$ and, by 2-goodness, $f''(m) > 0$.
Thus, the intersection of the convex hull with the second axis
consists of points with non-negative second coordinate.
From the proof of Lemma 3.5,
$\Sigma_1(u) = 0$ now implies $\Sigma_2(u) > 0$;
thus, the fold points have concavity upwards
in the restriction of $F$ to each fibre.
From the behaviour of $f$ at infinity,
$F$ has exactly one fold point per fibre.
The rest of the argument is similar to the proof of Theorem 4.2.
\qed

\goodbreak

{\nobf Corollary 4.4: }
{\sl
\meti{(a)}{Let $f:\RR \to \RR$ be a 2- and 3-good polynomial
of even degree and positive leading coefficient.
If $f''$ assumes both signs
then the operator $F$ has a cusp and there are
points in the image of $F$ with four pre-images.}
\meti{(b)}{There are non-convex smooth proper functions $f: \RR \to \RR$
for which $F$ is a global fold.}}

\goodbreak

{\nobf Proof: }
To prove (a), notice that,
by hypothesis, there is a point in the lower half-plane ($f'' < 0$).
Also, since $f$ is a polynomial, the argument of $\hat\gamma_2(t)$
approaches 0 from above (resp., $\pi$ from below)
when $t$ goes to $+\infty$ (resp., $-\infty$).
This suffices to show that 0
is in the interior of the convex hull of the image of $\hat\gamma_2(t)$.

As for (b), we deviate from the previous argument
by considering functions for which $f'$ and $f''$ are comparable for large $x$,
such as $\cosh x$.
More precisely, let $f$ be such that $f''$ 
coincides with $\cosh t$ outside a small interval of large positive numbers
in which we subtract from $\cosh t$ a narrow positive bump---%
if the bump is sufficiently high and narrow,
$f''$ changes sign, $f'(x) > 0$ for all positive $x$
and 0 is not in the convex hull of the image of $\hat\gamma_2$.
\qed

We can also characterize global folds in a purely local way:

{\nobf Theorem 4.5: }
{\sl If $f: \SS^1 \times \RR \to \RR$ is proper,
0 is a regular value of $\Sigma_1$ and all singularities
of $F$ are folds, then $F$ is a global fold.}

{\nobf Proof:}
For the operators $F$ being considered,
there is a sign associated to each fold:
it is the sign of $\Sigma_2$, which can also be interpreted
as saying whether the concavity of the restriction of $F$
to fibres points up or down.
This splits the set of folds into two open subsets.
Since by hypothesis $S_2 = \emptyset$ and by Corollary 1.9
$S_1$ is connected it follows that one of these sets is empty.
In other words, all folds are concave upwards, say,
and the restriction of $F$ to any fibre thus has at most
one critical point. The result follows by properness
and juxtaposition of fibres as in Theorems 4.2 and 4.3.
\qed

\bigbreak

{\noindent\bigbf 5. Global cusps}

\nobreak\smallskip\nobreak

The results in the previous section
describe the global geometry of $F$
when $f$ satisfies $D_2f > 0$ or $D_2^2f > 0$.
In the autonomous case,
we are able to handle another kind of nonlinearity.

{\nobf Theorem 5.1: }
{\sl Let $f: \RR \to \RR$ be a proper function such that
\item{\rm (a)}{$f'''(x) \ge 0$,}
\item{\rm (b)}{$f'''(x)$ has isolated roots,}
\item{\rm (c)}{$f'(x)$ assumes both signs.}

Then $F$ is a topological global cusp;
in the $H$ case, $F$ is a smooth global cusp.}

A similar result holds if (a) is replaced by $f'''(x) \le 0$.

The scenarios of the previous section are simple enough to allow
for rather explicit global changes of variable to normal form.
This is partly due to the fact that restrictions of $F$
to arbitrary fibres have similar behaviours:
in a sense, we may consider fibres individually.
For the operators in the statement of Theorem 5.1,
instead, such restrictions vary according to the fibre,
as illustrated in Figure 5.1:
we must therefore treat them collectively.
In the process, explicitness is lost:
from the theorem (and its proof), in the H case
the domain and image are foliated by smooth surfaces,
diffeomorphic to $\RR^2$, which are in turn foliated by fibres.
More, $F$ takes surfaces to surfaces
and, on each surface, $F$ is a global cusp.
Still, we have no idea how to exhibit such foliations:
they are shown to exist
by topological methods in Hilbert manifolds ([Ku], [MST]).
For Banach spaces, we use an additional existential
argument depending on the homeomorphism of all infinite dimensional separable
Banach spaces ([Ka]).

An operator $G: H^1 \to H^0$ is a {\sl smooth global cusp}
if there exist diffeomorphisms
$\Xi_1: H^1 \to H \times \RR^2$ and 
$\Xi_0: H^0 \to H \times \RR^2$ such that
$(\Xi_0 \circ G \circ \Xi_1^{-1})(Z, x, y) = (Z, x, y^3 + xy)$,
for all $(Z, x, y) \in H \times \RR^2$,
where $H$ is a separable infinite dimensional Hilbert space.
Similarly, we call $G: C^1 \to C^0$ a {\sl topological} global cusp
if there exist homeomorphisms $\Xi_i: C^i \to H \times \RR^2$ as above.

The following lemma is an exercise in integration by parts.

{\nobf Lemma 5.2: }
{\sl Let $g: \RR \to \RR$ be a smooth function.
Given $a < b$, we have
$$(b-a)(g'(a) + g'(b)) - 2(g(b) - g(a)) =
- \int_a^b (t-a)(t-b) g'''(t) dt.$$
}

An additional topological ingredient is an infinite-dimensional version
of the corollary to Theorem 1 in [S].
The proof is similar for Hilbert spaces and makes use of results in [BH]
for Banach spaces (see [MST] for additional information).

{\nobf Lemma 5.3: }
{\sl 
\meti{\rm (a)}
Given a contractible connected smooth submanifold
$H'$  of codimension 1 of a separable Hilbert space $H$ of infinite
dimension, there is a diffeomorphism of $H$ to itself taking $H'$
to a closed subspace of codimension 1.
\meti{\rm (b)}
Let $B'$ be a closed subset
of a separable infinite-dimensional Banach space $B$.
Assume that $B'$ is connected, contractible,
and a bicollared topological submanifold of codimension 1 in $B$.
Then there is a homeomorphism from $B$ to itself
taking $B'$ to a closed subspace of codimension 1.
}

Actually, a contractible connected closed bicollared topological submanifold of
codimension 1 always splits the ambient space in (exactly) two components.

Finally, we make use of a canonical construction
to bring planar global cusps to normal form.
A sketch of proof is given at the end of this section.

{\nobf Lemma 5.4: }
{\sl
Let $\cal Z$ be a topological space of parameters.
Let $G: {\cal Z} \times \RR^2 \to {\cal Z} \times \RR^2$
be a continuous function of the form
$$G(Z,x,y) = (Z,x,g_Z(x,y))$$
with the following properties:
\itembu{For all $Z \in {\cal Z}$, $g_Z(0,0) = 0$.}
\itembu{For all $Z \in {\cal Z}$ and $x \in \RR$,
$\lim_{y \to \pm \infty} g_Z(x,y) = \pm \infty$.}
\itembu{For all $Z \in {\cal Z}$ and $x \le 0$ the function
$y \mapsto g_Z(x,y)$ is strictly increasing.}
\itembu{There exist continuous functions
$m, M: {\cal Z} \times [0, +\infty) \to \RR$
with $m(Z,0) = M(Z,0) = 0$ such that,
for all $Z \in {\cal Z}$ and $x > 0$ the function
$y \mapsto g_Z(x,y)$ is strictly increasing in $(-\infty, m(Z,x)]$
and $[M(Z,x), +\infty)$ but strictly decreasing in $[m(Z,x),M(Z,x)]$.}

\meti{\rm (a)}
There exist homeomorphisms $W_d$ and $W_i$ of ${\cal Z} \times \RR^2$
keeping the $Z$ and $x$ coordinates fixed such that
$$(W_i \circ G \circ W_d^{-1})(Z,x,y) = (Z,x,y^3 + xy).$$
\meti{\rm (b)}
If ${\cal Z}$ is a smooth Hilbert manifold
and $G$ is smooth with 
$$D_2g_Z(0,0) = D_2^2g_Z(0,0) = 0,\quad D_2^3g_Z(0,0) > 0,
\quad D_1D_2g_Z(0,0) < 0$$
then $W_d$ and $W_i$ can be taken to be diffeomorphisms.
}

{\nobf Proof of Theorem 5.1: }

We first classify the singularities, next we study
the behaviour of the restriction of $F$ to fibres
and finally obtain global results by juxtaposing fibres.

{\sl All critical points of $F$ are folds or cusps.}
Clearly, $\Sigma_3$ is always positive and we only have
to check transversality when $\Sigma_1 = \Sigma_2 = 0$,
i.e., we have to verify that $D\Sigma_1 \ne 0$ at such points.
Since $\Sigma_1 = \ints f'(u(t)) dt$,
$$D\Sigma_1(u) \cdot v = \ints f''(u(t)) v(t) dt;$$
we show that the function $f''(u(t))$ is not identically zero.
First, there is a unique $x_0$ for which $f''(x_0) = 0$;
from (c), we have $f'(x_0) < 0$.
Thus, the function $f''(u(t))$ is identically zero
only if $u(t) = x_0$ for all $t$ but then
$\ints f'(x_t) dt < 0$ and $u \notin S_1$.
\qede

We now consider $F$ restricted to fibres.
From Theorem 1.4, such restrictions take $\pm \infty$ to $\pm \infty$.
Also, regular and fold points of $F$ are regular or fold points
of the restriction.
Furthermore, the restrictions are locally increasing
at cusp points.
Indeed, $\Sigma_3 > 0$ and, at cusp points,
$\Sigma_3$ is a positive multiple of $\wz^3\Phi$,
the third derivative of the restriction.

{\sl The operator $F$ has at most two critical points per fibre.}
Let $u_1 < u_2$ be two critical points of $F$ in the same fibre:
$$\ints f'(u_1(t)) dt = \ints f'(u_2(t)) dt = 0\eqno{(\ast)}$$
and
$$u_1'(t) + f( u_1(t)) = u_2'(t) + f( u_2(t)) + C\eqno{(\dag)}$$
for some constant $C$. We first show that $C > 0$.
From Lemma 5.2,
$$f'(u_2(t)) + f'(u_1(t)) -
{f(u_2(t)) - f(u_1(t)) \over u_2(t) - u_1(t)} > 0$$
and, integrating and making use of $(\ast)$,
$$\ints {f(u_2(t)) - f(u_1(t)) \over u_2(t) - u_1(t)} dt < 0.$$
From $(\dag)$,
$$ {(u_2 - u_1)'(t) \over (u_2 - u_1)(t)} +
{f(u_2(t)) - f(u_1(t)) \over u_2(t) - u_1(t)} +
{C \over u_2(t) - u_1(t)} = 0$$
and, integrating, we obtain
$$C \ints {dt \over u_2(t) - u_1(t)} > 0$$
implying $C>0$.
This means that the restriction of $F$ to fibres,
if further restricted to critical points, is decreasing:
$u_1 < u_2$ implies $F(u_1) > F(u_2)$.
It follows that, if there are at least three critical points
on a fibre, $F$ is decreasing near the second one:
from the previous paragraph,
this second critical point can thus be neither a fold nor a cusp.
\qede

At this point we already know that regular values have one or three
pre-images, images of folds have two pre-images 
and images of cusps have a single pre-image (as in [CD] and [CT]).
But we have more: $F$ restricted to a fibre
is topologically equivalent to one of the three graphs on Figure 5.1.
Fibres containing a cusp, on which $F$ behaves as depicted in (b),
split the space of fibres into two subsets,
on which restrictions of $F$ behave as in either (a) 
(no critical points) or (c) (two folds).
Let $\Fx = \tilde B^1$ be, as before, the space of fibres:
we have the natural partition $\Fx = \Fa \cup \Fb \cup \Fc$
into sets of fibres of types (a), (b) and (c), respectively.

\midinsert
\line{\hfill%
\vbox{\hsize=4cm \psfig{file=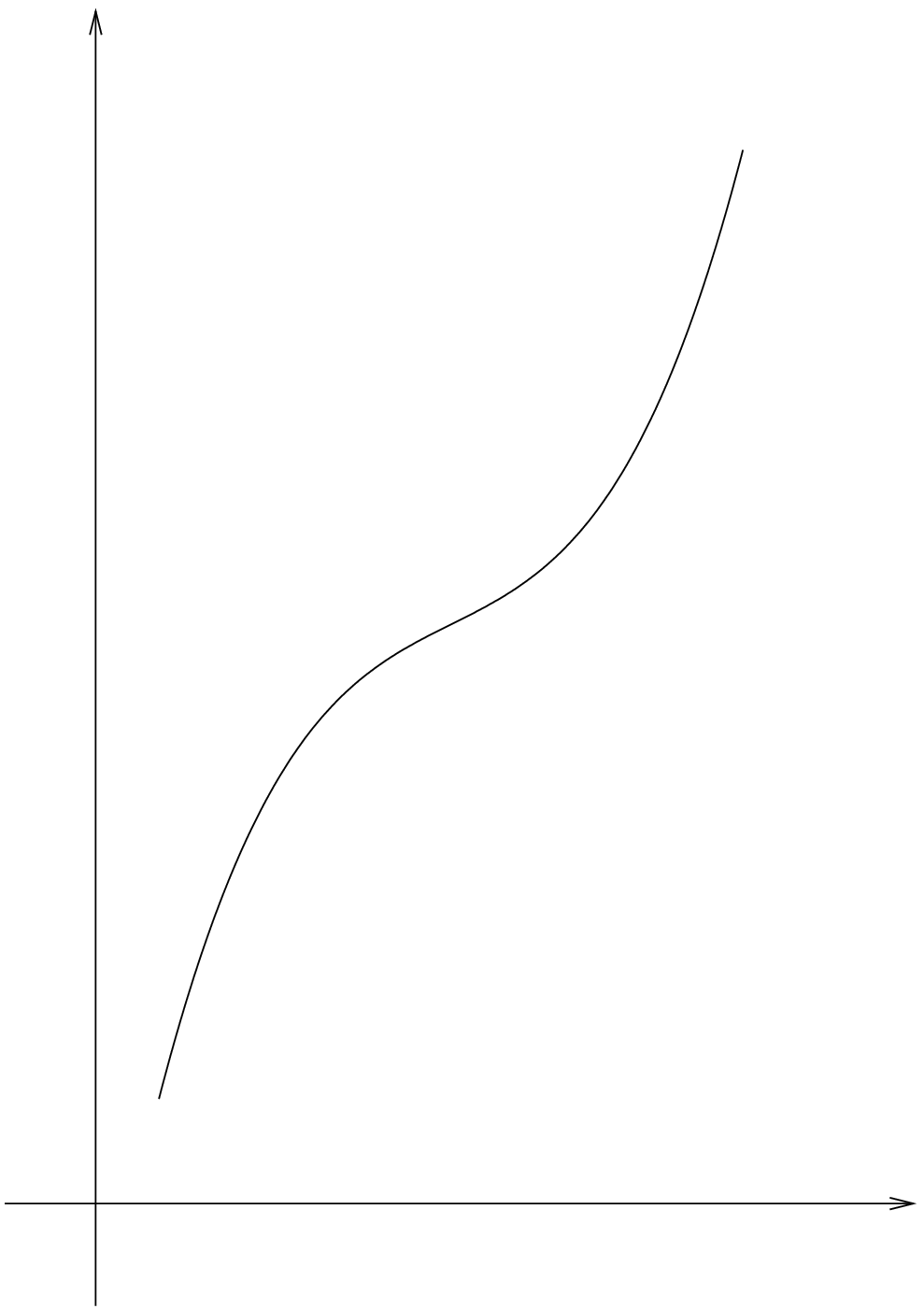,width=4cm} \smallskip\centerline{(a)}}%
\quad%
\vbox{\hsize=4cm \psfig{file=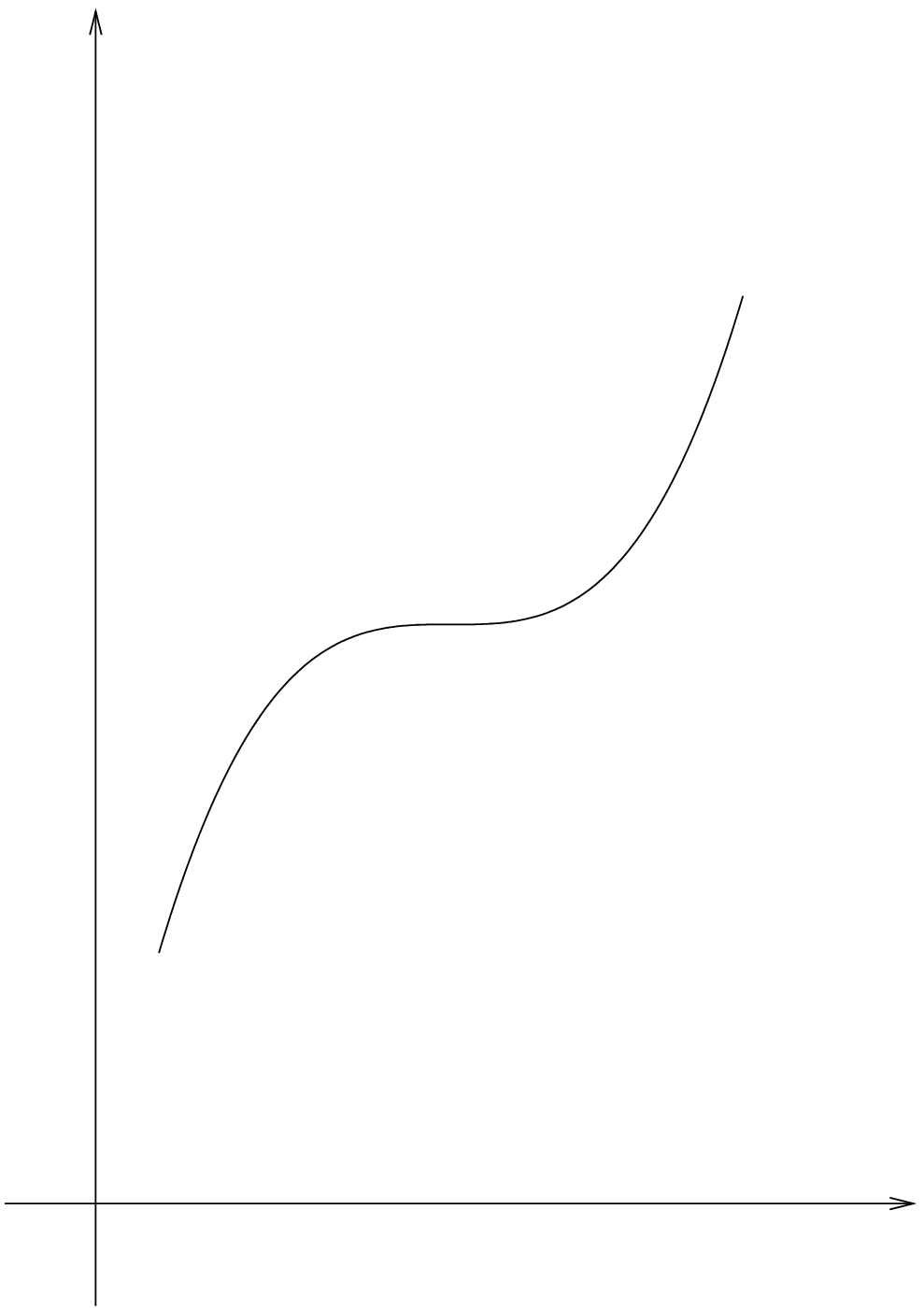,width=4cm} \smallskip\centerline{(b)}}%
\quad%
\vbox{\hsize=4cm \psfig{file=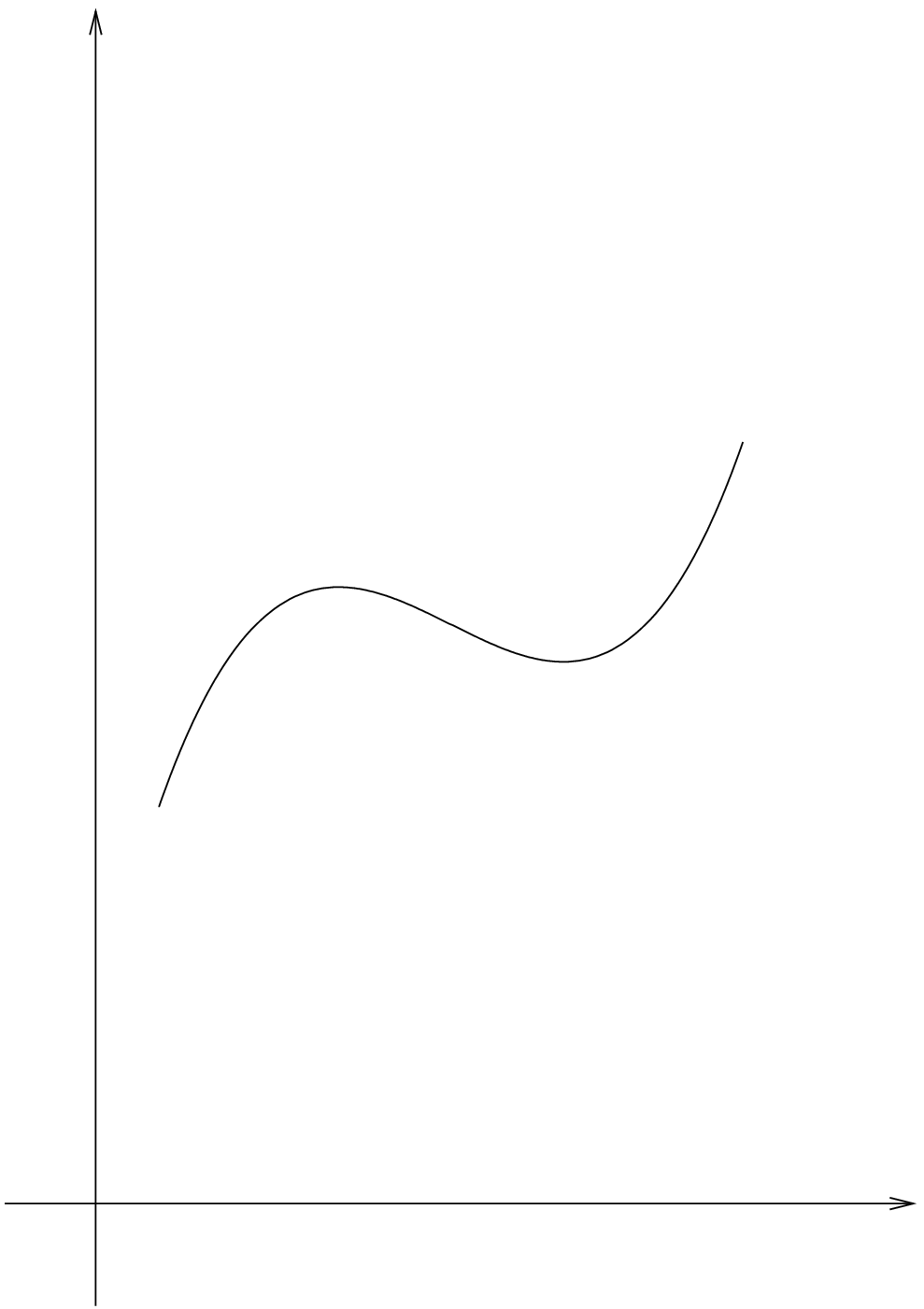,width=4cm} \smallskip\centerline{(c)}}%
\hfill}
\smallskip
\centerline{\eightrm Figure 5.1}
\endinsert

{\sl The sets $\Fa$ and $\Fc$ are open in $\Fx$; $\Fb$ is closed.}
Let $u_0 \in \Fc$;
there exist elements $u_1 < u_2$ in the $u_0$-fibre with $F(u_1) > F(u_2)$.
The hyperplanes parallel to $\tilde B^1$ passing
through $u_1$ and $u_2$ transversally intersect fibres.
In particular, fibres sufficiently close to the $u_0$-fibre
contain point $u'_1 < u'_2$ in these hyperplanes 
for which $F(u'_1) > F(u'_2)$, proving the openness of $\Fc$.

Assuming by contradiction that $\Fa$ is not open,
let $u_n$ be a sequence of critical points
whose corresponding fibres converge to the fibre $u_\infty \in \Fa$.
If the averages of $u_n$ are bounded,
we may assume by compactness that these averages converge;
this, however, implies that the sequence $u_n$ itself
converges to the element of the $u_\infty$-fibre with limiting average.
Since this limit is clearly a critical point of $F$
we have the contradiction in this case
and may assume from now on that the averages of $u_n$
tend monotonically to $+\infty$.

Let $M$ be such that $x > M$ implies $f'(x) > 0$.
Notice in particular that for all $u \in S_1$
there is a $t \in \SS^1$ such that $f'(u(t)) = 0$
and therefore $u(t) < M$.
For each $n$, let $t_n \in \SS^1$ be such that $u_n(t_n) < M$.
From the compactness of $\SS^1$, we may assume that the sequence
$t_n$ converges to, say, $t_\infty$. 
Let $u_{nm}$ ($m < n$) be the element in the $u_n$-fibre 
with average $\overline u_m$; clearly, $u_{nm} < u_n$.
Define similarly $u_{\infty m}$:
the sequence $u_{nm}$ (for fixed $m$) tends to $u_{\infty m}$.
Since $u_{nm}(t_n) < M$ for all $n$,
$u_{\infty m}(t_\infty) \le M$ for all $m$,
in contradiction with Lemma 1.3.
\qede

{\sl The set $S_2$ of cusps of $F$ is a smooth contractible 
submanifold of codimension 2 of $B^1$.}
From Lemma 3.1, there is a diffeomorphism of $B^1$ to itself
taking the sets $S_1$ and $S_2$ to $S_1 = \hat S_1$ and $\hat S_2$,
respectively and it suffices to show that $\hat S_2$ is contractible.
We may thus use Corollary 3.3: we check that $\hat F$ is 2-regular,
i.e., that 0 is a regular value for both $\hat \Sigma_1$
and $(\hat \Sigma_1, \hat \Sigma_2)$.
We have already seen that $D\Sigma_1 = D(\hat \Sigma_1)$
is never zero in the critical set.
Also,
$$D(\hat \Sigma_1, \hat \Sigma_2)(u) \cdot v =
\left(\ints f''(u(t)) v(t) dt, \ints f'''(u(t)) v(t) dt\right) $$
and we are left with showing the linear independence of the functions
$f''(u(t))$ and $f'''(u(t))$ for $u$ satisfying
$\ints f'(u(t)) dt = \ints  f''(u(t)) dt = 0$.
The first function is non-zero but of average zero
and the second is strictly positive.
\qede

In particular, from Corollary 3.3, there is a homeomorphism
of $B^1$ to itself taking $S_1$ and $S_2$ to nested subspaces
of codimensions 1 and 2.
This homeomorphism, however, does not respect fibres:
we now show how to do better.
It is convenient from this point on to work in adapted coordinates,
i.e., to consider $\Fbd: B^0 \to B^0$, its critical set $\Sbd_1 = \Psi(S_1)$
and set of cusps $\Sbd_2 = \Psi(S_2)$.
Recall that $\Psi: B^1 \to B^0$ takes fibres to vertical lines
and that $\Fbd$ is in adapted coordinates:
in other words, vertical lines are fibres for $\Fbd$.
Similarly, $\tilde B^0 = \Fabd \cup \Fbbd \cup \Fcbd$,
the disjoint images under $\Psi$ of $\Fa$, $\Fb$ and $\Fc$.

In the next claims, we construct a number of auxiliary
changes of variable, leading eventually to the global normal form for the cusp.
In the H case, all constructions are smooth.
In the C case, however, we make use of homeomorphisms
and folds and cusps have to be interpreted topologically.
Figure 5.2 may be helpful.

\midinsert
\centerline{\psfig{width=4.5in,file=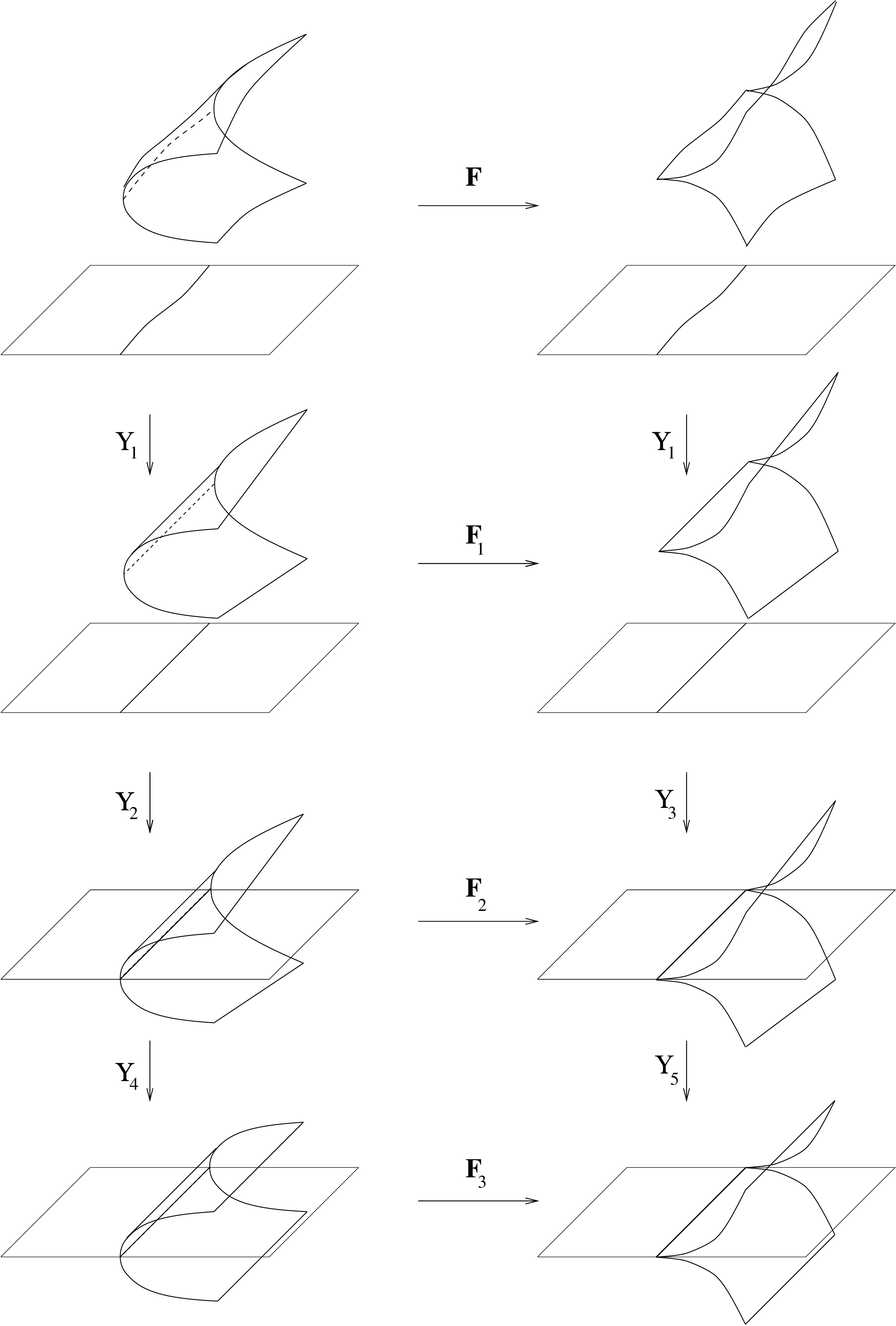}}
\smallskip
\centerline{\eightrm Figure 5.2}
\endinsert

{\sl There is a homeomorphism $\Upsilon_1$ of $B^0$ to itself
taking vertical lines to vertical lines
and $\Fbbd$ to a subspace $\Cx$ of codimension 1 of $\tilde B^0$.
In the H case, this homeomorphism can be taken to be a diffeomorphism.}
From the local normal form for cusps,
vertical lines intersect $\Sbd_2$ transversally:
the vertical projection is then a natural diffeomorphism
between $\Sbd_2$ and $\Fbbd$.
In particular, $\Fbbd \subset \tilde B^0$ is a contractible submanifold
of codimension 1 for which Lemma 5.3 applies:
smoothness guarantees the existence of local tubular neighbourhoods
which can be consistently glued because the complement of $\Fbbd$
has two connected components.
Thus, there is a homeomorphism
$\Upsilon_0: \tilde B^0 \to \tilde B^0$
taking $\Fbbd$ to a closed subspace $\Cx$ of codimension 1.
Define $\Upsilon_1: B^0 \to B^0$ as the only extension of $\Upsilon_0$
respecting vertical lines and horizontal hyperplanes%
---$\Upsilon_1$ is clearly a homeomorphism.
\qede

Notice that the conjugation $\Fbd_1 =
\Upsilon_1 \circ \Fbd \circ \Upsilon_1^{-1}$
is still in adapted coordinates, i.e., vertical lines are invariant.
Also, the vertical projection of $\Upsilon_1(\Sbd_2)$
and $\Upsilon_1(\Fbd(\Sbd_2))$ are both equal to $\Cx$;
$\Upsilon_1(\Sbd_2)$ is the set of cusps (or, in the C case,
topological cusps) of $\Fbd_1$ and $\Upsilon_1(\Fbd(\Sbd_2))$ is its image.

Identify
$$B^0 = \Cx \oplus \langle r \rangle \oplus \langle 1 \rangle,$$
for $r \in \tilde B^0$ not in $\Cx$ and now
write a typical element of $B^0$ as $(Z,x,y) \in \Cx \times \RR^2$
making use of the natural projections.
Notice that the sign of $x$ determines the type of the $(Z,x,0)$-fibre:
without loss, the cases $x < 0$, $x = 0$ and $x > 0$
are set to correspond to fibres of types (a), (b) and (c), respectively.

{\sl There are homeomorphisms $\Upsilon_2, \Upsilon_3$ of $B^0$ to itself
keeping each vertical line invariant and such that
$\Upsilon_2(\Upsilon_1(\Sbd_2)) = 
\Upsilon_3(\Upsilon_1(\Fbd(\Sbd_2))) = \Cx$.
In the H case, these maps are diffeomorphisms.}
Set $\Upsilon_2(Z,x,y) = (Z,x,y-y')$
where $y'$ is the only real number such that
$\Upsilon_1^{-1}(Z,0,y') \in \Sbd_2$.
Similarly, set $\Upsilon_3(Z,x,y) = (Z,x,y-y'')$
where $y''$ satisfies
$\Upsilon_1^{-1}(Z,0,y'') \in \Fbd(\Sbd_2)$.
\qede

The composition $\Fbd_2 = \Upsilon_3 \circ \Fbd_1 \circ \Upsilon_2^{-1}$
is almost in normal form: for each $Z$, $\Fbd_2$ restricted to the $Z$-plane
is a global 2-dimensional cusp.
We are now ready to apply Lemma 5.4 to get two further changes
of variable $\Upsilon_4$ and $\Upsilon_5$ such that
$\Fbd_3 = \Upsilon_5 \circ \Fbd_2 \circ \Upsilon_4^{-1}$
is the desired normal form
$$\Fbd_3(Z,x,y) = (Z,x,y^3 + xy).$$
\qed

{\nobf Sketch of proof of Lemma 5.4: }
Item (a) is straightforward.
As to item (b), we begin by invoking Whitney's construction ([W])
which obtains, for any given $Z$,
diffeomorphisms $\upsilon_{Z,6}$ and $\upsilon_{Z,7}$
of neighbourhoods
of the origin taking the restriction $G_Z(x,y) = (x,g_Z(x,y))$
to the normal form $k(x,y) = (x, y^3 + xy)$,
i.e., $k = \upsilon_{Z,7} \circ G_Z \circ \upsilon_{Z,6}^{-1}$.
The diffeomorphisms $\upsilon_{Z,6}$ and $\upsilon_{Z,7}$
so constructed are of the form
$(x,y) \mapsto (x,y')$ and preserve orientation.

Actually (and here we omit the verification),
the construction allows $\upsilon_{Z,6}$, $\upsilon_{Z,7}$
and the size of the neighbourhoods to be chosen smoothly
as functions of the parameter $Z$.
More exactly, we have diffeomorphisms $\Upsilon_{T,6}$ and $\Upsilon_{T,7}$
defined on tubular neighbourhoods of ${\cal Z} \times (0,0)$
taking $G$ to the normal form 
$K(Z,x,y) = (Z,x,y^3 + xy)$ 
near ${\cal Z} \times (0,0)$,
i.e., $K = \Upsilon_{T,7} \circ G \circ \Upsilon_{T,6}^{-1}$
whenever the right hand side is defined.
Now extend $\Upsilon_{T,6}$ and $\Upsilon_{T,7}$
to diffeomorphisms $\Upsilon_6$ and $\Upsilon_7$
from ${\cal Z} \times \RR^2$ to itself
of the form $(Z,x,y) \mapsto (Z,x,y')$
which coincide with the identity outside a
tubular neighbourhood of ${\cal Z} \times (0,0)$---%
notice that this extension is just a one-dimensional problem,
parametrized by $(Z,x)$.
The composition $G_1 = \Upsilon_{7} \circ G \circ \Upsilon_{6}^{-1}$
satisfies all the original conditions in Lemma 5.4
and coincides with the normal form $K$
in a tubular neighbourhood of ${\cal Z} \times (0,0)$.
The hard part of bringing the cusps into normal form being done,
we give explicit instructions to take $G_1$ to normal form.

Let $S_{G_1}$ and $S_K$ be the critical sets of $G_1$ and $K$.
Their images $G_1(S_{G_1})$ and $K(S_K)$
have two points in each vertical line $(Z,x,\cdot)$, $x > 0$.
Construct $\Upsilon_9$ by juxtaposing 1-dimensional maps
to be a diffeomorphism of ${\cal Z} \times (0,0)$
satisfying
\item{$\cdot$}$\Upsilon_9(Z,x,y) = (Z,x,y')$,
\item{$\cdot$}for each $Z$ there is a positive $x_Z$
such that $\Upsilon_9(Z,x,y) = (Z,x,y)$ for $x < x_Z$,
\item{$\cdot$}$\Upsilon_9(G_1(S_{G_1})) = K(S_K)$.

Let $G_2 = \Upsilon_9 \circ G_1$; notice that $G_2(S_{G_2}) = K(S_K)$,
where $S_{G_2}$ is the critical set of $G_2$.
There is now a unique diffeomorphism $\Upsilon_8$
satisfying $G_2 \circ \Upsilon_8^{-1} = K$.
Indeed, if $x \le 0$, let
$\Upsilon_8(Z,x,y)$ be the only point $(Z,x,y')$
for which $G_2(Z,x,y) = K(Z,x,y')$.
For positive $x$, given $(Z,x,y)$ there are at most three points
satisfying $G_2(Z,x,y) = K(Z,x,y')$:
we define $\Upsilon_8(Z,x,y)$
to be the only point $(Z,x,y')$ which is in the same position 
with respect to the two critical points of $K$
in the vertical line $(Z,x,\cdot)$
as $(Z,x,y)$ is with respect to the two critical points of $G_2$
in the same vertical line.
Clearly, $\Upsilon_8$ is a homeomorphism:
we check smoothness of $\Upsilon_8$ and its inverse.
At regular points, this is the inverse function theorem.
At fold points, one may use the square root trick in the proof of Theorem 4.2,
but we omit the details.
Finally, smoothness at cusps is guaranteed from the simple fact
that $\Upsilon_8$ turns out to be the identity
near ${\cal Z} \times (0,0)$.
\qed

As a corollary, we obtain a global cusp form
for the Cafagna-Donati operator ([CD], [CT]):

{\nobf Corollary 5.5: }
{\sl Let $f(x) = ax + bx^2 + cx^{2k+1}$ where
$k$ is a positive integer, $a \ge 0$, $a^2 + b^2 > 0$
and $c < 0$. Then the operator $F: H^1 \to H^0$
is a smooth global cusp and $F: C^1 \to C^0$
is a topological global cusp.}

\bigbreak

{\noindent\bigbf 6. A numerical counter-example}

\nobreak\smallskip\nobreak

It is of course tempting to speculate about
the possible consequences of $D_2^4f > 0$:
does this condition at least
guarantee that points have at most four pre-images?
In [L], Lins Neto shows that,
if $f(x,t)$ is a polynomial of degree four in $x$
with coefficients depending on $t$
and positive highest degree coefficient,
then the number of solutions may be arbitrarily large.

In this section, 
we obtain a polynomial $f$ of degree 4 and
a smooth periodic $u_b$ such that $u_b$
is a Morin singularity of order 4 (a {\sl butterfly}).
From Morin's normal form,
some points $v$ near $F(u_b)$ have five regular pre-images
close to $u_b$.
Since the degree of $F$ is zero, there is yet a sixth pre-image
and we thus obtain a smooth periodic function $v$
for which the equation
$$u' + f(u) = v(t), u(0) = u(1),$$
has six solutions.

We briefly describe the numerical procedure
used in the identification of $u_b$.
Without loss of generality, $f(x) = x^4 - b x^2 + c x$
and we try to find a butterfly $u_b$ of the form
$$u_b(t) = a_0 + a_1 \cos(t) +
a_2 \cos(2t) + b_2 \sin(2t) + \cdots
+ b_4 \sin(2t).$$
We now write the four scalar equations $\Sigma_i(u_b) = 0$, $i = 1, \ldots, 4$,
in terms of the ten parameters $b, c, a_0, \ldots, b_4$
and search for a zero with a Newton-like method
with pseudo-inverses [AG].
Actually, for appropriate $b$ and $c$, four extra parameters
should be enough to locate a butterfly, but the numerical analysis
becomes more robust with additional parameters.
One example is $b = 4$, $c = -0.3$,
$a_0 = -0.01173378$, $a_1 = -0.8836063$, 
$a_2 = 0.2428734$, $b_2 = -0.6855379$, 
$a_3 = 0.4465347$, $b_3 = 0.1853376$, 
$a_4 = -0.01881213$ and $b_4 = 0.2105862$.
The Newton method itself checks for the surjectivity
of (the restriction of) the derivative of $(\Sigma_1,\ldots,\Sigma_4)$
and the program also verifies that $\Sigma_5(u_b) \ne 0$.

\midinsert
\centerline{\psfig{file=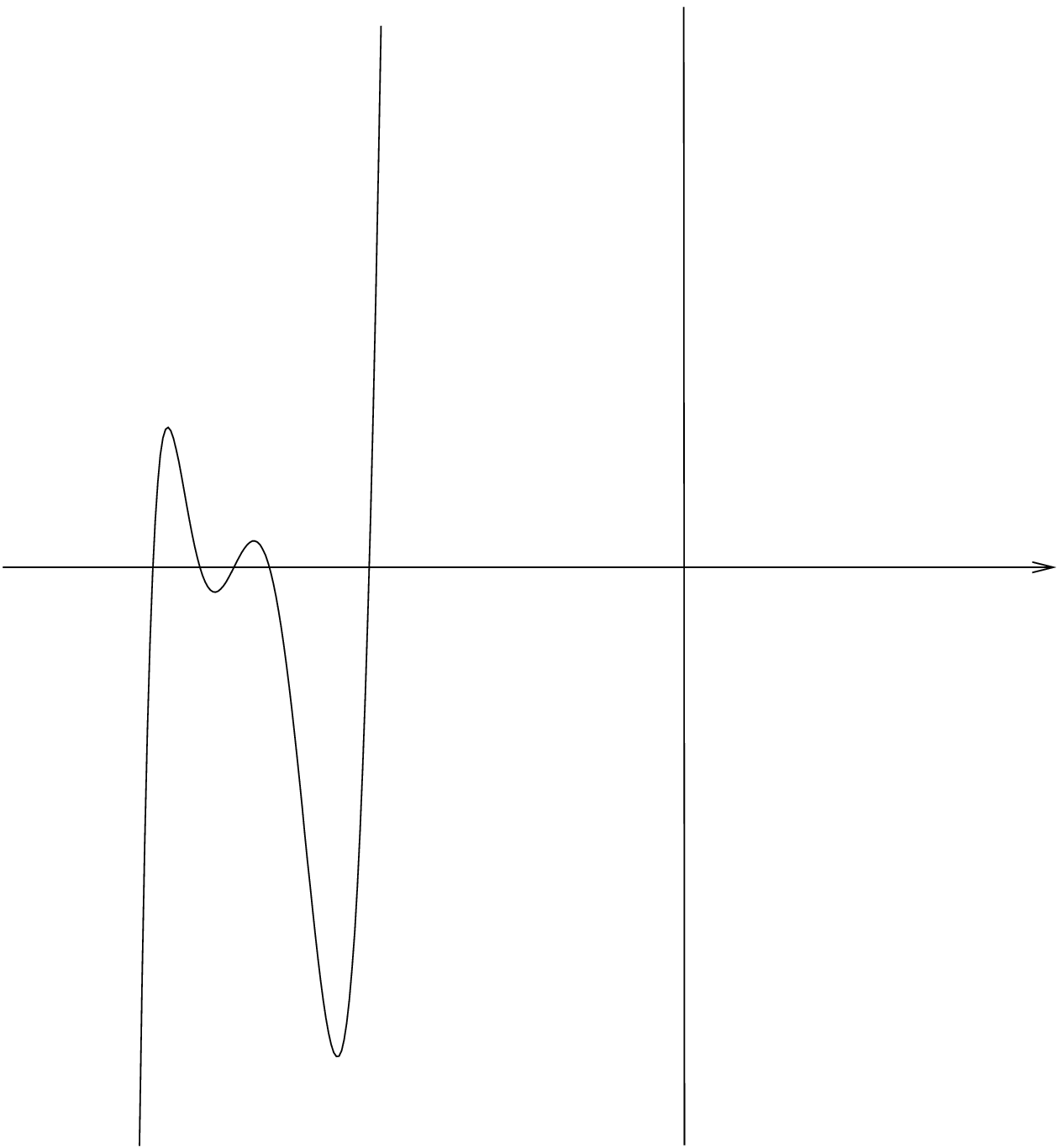,width=3in}}
\smallskip
\centerline{\eightrm Figure 6.1}
\endinsert

Again by Newton's method, we try to solve
$$(\Sigma_1,\ldots,\Sigma_4)(u_1) = (-0.0000005, 0, 0.00008, 0),$$
where the non-zero constants on the right hand side
were adjusted somewhat empirically---in a nutshell,
we are trying to perturb the polynomial $x^5$
to get five distinct real roots, which can be accomplished
by adding small multiples of $x^3$ and $x$.
The parameters for $u_1$ are
$$\matrix{
&a_0 = -0.011367708203969,&
a_1 = -0.883600656945802, \cr
a_2 = 0.243308077825844,&
a_3 = 0.446085678376277,&
a_4 = -0.018458472190807,\cr
b_2 = -0.685621717642052, &
b_3 = 0.185481811055651, &
b_4 = 0.210509692732880.\cr
}$$
In Figure 6.1 we plot $\rho_v(x) - x$, where $v = u_1' + f(u_1) = F(u_1)$,
so that roots of this auxiliary function
correspond to periodic solutions of $u' + f(u) = v$.
This graph was obtained by solving the differential equation
with a Runge-Kutta method of order 4 for initial conditions ranging from
$-0.4$ to $0.4$.
The vertical scale is stretched by a factor of $2.5 \cdot 10^6$
and the time step for the method had to be taken as $2 \cdot 10^{-4}$.
Notice the clustering of the first five roots,
stemming from the butterfly:
actually it is this clustering and the quintic behaviour of the butterfly
which account for the need of a huge vertical stretching factor.
Another consequence of the quintic behaviour is the great
sensitivity of the coefficients:
for instance, a change of $10^{-6}$ in $a_0$ destroys four of the six solutions.
Still, this final direct check is far easier (and more reliable)
than the process of obtaining the coefficients for the example.

\bigbreak

{\noindent\bigbf References: }

\nobreak\medskip\nobreak

\parindent=40pt

\item{[AG]}{ Allgower, E. L. and Georg, K.,
{\sl Numerical continuation methods: an introduction},
Springer-Verlag, New York, 1991.}

\item{[AP]}{ Ambrosetti, A. and Prodi, G.,
{\sl On the inversion of some differentiable maps between Banach spaces
with singularities},
Ann. Mat. Pura Appl. 93 (1972), 231-246.}

\item{[BC]}{ Berger, M. S. and Church, P. T.,
{\sl Complete integrability and a perturbation of a non-linear
Dirichlet problem (I)},
Indiana Univ. Math. J., 28 (1979), 935-952.}

\item{[BCT]}{ Berger, M. S., Church, P. T. and Timourian, J. G.,
{\sl Folds and cusps in Banach spaces, with applications to
nonlinear partial differential equations (I)},
Indiana Univ. Math. J., 34 (1985), 1-19.}

% \item{[BH]}{ Burghelea, D. and Henderson, D.,
% {\sl Smoothings and homeomorphisms for Hilbert manifolds},
% Bull. Amer. Math. Soc., 76 (1970), 1261--1265.}

\item{[CD]}{ Cafagna, V. and Donati, F.,
{\sl Un r\'esult global de multiplicit\'e pour un probl\`eme diff\'erentiel
non lin\'eaire du premier ordre},
C. R. Acad. Sc. Paris, 300 (1985), 523--526.}

\item{[CDT]}{ Church, P. T., Dancer, E. N. and Timourian, J. G.,
{\sl The structure of a nonlinear elliptic operator},
Trans. Amer. Math. Soc., 338, no. 1 (1993), 1--42.}

\item{[CT]}{ Church, P. T. and Timourian, J. G.,
{\sl Global cusp maps in differential and integral equations},
Nonlinear Analysis, 20, no. 11 (1993), 1319--1343.}

\item{[G]}{ Gr\"unbaum, B.,
{\sl Convex polytopes},
Interscience, London, 1967.}

\item{[Ka]}{ Kadec, M. I.,
{\sl A proof of the topological equivalence
of all separable infinite-dimensional Banach spaces},
Funkcional Anal. i Prilozen, 1 (1967), 61--70.}

\item{[Ku]}{ Kuiper, N. H.,
{\sl Vari\'et\'es Hilbertiennes --- Aspects G\'eom\'etriques},
Publications du S\'eminaire de Math\'ematiques Sup\'erieures, No. 38,
Les Presses de l'Universit\'e de Montreal, 1971.}

\item{[L]}{ Lins Neto, A.,
{\sl On the number of solutions of the equation
${dx \over dt} = \sum_{j=0}^n{a_j(t) x^j}$, $0 \le t \le 1$,
for which $x(0) = x(1)$},
Inventiones Math., 59 (1980), 67--76.}

%\item{[L]}{ Lang, S., {\sl Differentiable Manifolds},
%Addison-Wesley, Reading, MA, 1972.}

\item{[M]}{ Morin, B.,
{\sl Formes canoniques de singularit\'es d'une application diff\'erentiable},
C. R. Acad. Sc. Paris, 260 (1965), 5662--5665 and 6503--6506.}

\item{[McKS]}{ McKean, H. P., Scovel, J. C.,
{\sl Geometry of some simple nonlinear differential operators},
Ann. Scuola Norm. Sup. Pisa Cl. Sci. (4) 13 (1986), no. 2, 299--346.}

\item{[MST]}{ Malta, I., Saldanha, N. C. and Tomei, C.,
{\sl Regular level sets of averages of
Nemytski{\u\i} operators are contractible},
to appear in J. Func. Anal.}

% Steinitz, E., Bedingt konvergente Reihen und konvexe Systeme,
% J. reine angew. Math., ...

\item{[S]}{ Stallings, J. R.,
{\sl On Infinite Processes Leading to Differentiability
in the Complement of a Point}, 245--254,
in {\sl Differentiable and Combinatorial Topology,
a symposium in honor of Marston Morse},
ed. Cairns, S. S.,
Princeton University Press, Princeton, NJ, 1965.}

\item{[W]}{Whitney. H.,
{\sl On singularities of mappings of Euclidean spaces I.
Mappings of the plane into the plane},
Ann. Math., 62 (1955), 374--410}.

\bigskip\bigskip\bigbreak

{

\parindent=0pt
\parskip=0pt
\obeylines

Iaci Malta and Carlos Tomei, Departamento de Matem\'atica, PUC-Rio
R. Marqu\^es de S. Vicente 225, Rio de Janeiro, RJ 22453-900, Brazil

\medskip

Nicolau C. Saldanha, IMPA
Estr. D. Castorina 110, Rio de Janeiro, RJ 22460-320, Brazil

\medskip

malta@mat.puc-rio.br
nicolau@impa.br; http://www.impa.br/$\sim$nicolau/
tomei@mat.puc-rio.br

}

\bye